\documentclass[11pt,a4paper]{article}
\usepackage{graphicx}
\usepackage{amsmath}
\usepackage{amssymb}
\usepackage{amsthm}
\usepackage{cite}
\usepackage{subfigure}

\newcommand{\bm}{\boldsymbol}
\newcommand{\vect}[1]{\mathbf{#1}}
\newcommand{\mat}[1]{\mathbf{#1}}

\newcommand{\E}[1]{{\mathbb{E}}\{#1\}}

\newcommand{\etr}[1]{{\mathrm{etr}}\left\{#1\right\}}

\newcommand{\diag}[1]{{\mathrm{diag}}(#1)}
\renewcommand{\det}[1]{|#1|}
\newcommand{\Prob}[1]{{\mathbb{P}}\left(#1\right)}
\newcommand{\Proba}[2]{{\mathbb{P}}_{#1}\left(#2\right)}
\newcommand{\chol}[1]{\mathrm{chol}\left(#1\right)}

\newcommand{\ttilde}{\tilde{t}}
\newcommand{\xtilde}{\tilde{x}}
\newcommand{\Pfa}{P_{fa}}
\newcommand{\Pfabar}{\bar{P}_{fa}}
\newcommand{\Pd}{P_{d}}
\newcommand{\Pdbar}{\bar{P}_{d}}
\newcommand{\gammat}{\gamma_{t}}

\newcommand{\elast}{\vect{e}_{N}}
\newcommand{\vl}{\boldsymbol{\ell}}
\newcommand{\vn}{\vect{n}}

\newcommand{\vq}{\vect{q}}
\newcommand{\vx}{\vect{x}}
\newcommand{\vxbar}{\bar{\vx}}
\newcommand{\vxtilde}{\tilde{\vx}}

\newcommand{\vv}{\vect{v}}
\newcommand{\vvbar}{\bar{\vv}}

\newcommand{\vzero}{\vect{0}}

\newcommand{\F}{\mat{F}}
\newcommand{\Ft}{\F_{t}}
\newcommand{\G}{\mat{G}}
\newcommand{\Gt}{\G_{t}}
\newcommand{\I}{\mat{I}}
\newcommand{\eye}[1]{\I_{#1}}

\newcommand{\Q}{\mat{Q}}
\newcommand{\Qv}{\Q_{\mathrm{v}}}

\newcommand{\St}{\mat{S}_{t}}
\newcommand{\U}{\mat{U}}
\newcommand{\V}{\mat{V}}
\newcommand{\Vorth}{\V_{\perp}}
\newcommand{\W}{\mat{W}}
\newcommand{\Wt}{\W_{t}}
\newcommand{\X}{\mat{X}}
\newcommand{\Xbar}{\bar{\X}}
\newcommand{\Xt}{\X_{t}}

\newcommand{\Z}{\mat{Z}}
\newcommand{\Mzero}{\mat{0}}

\newcommand{\mLambda}{\bm{\Lambda}}

\newcommand{\mSigma}{\bm{\Sigma}}
\newcommand{\mSigmat}{\bm{\Sigma}_{t}}
\newcommand{\mOmega}{\bm{\Omega}}
\newcommand{\mPsi}{\bm{\Psi}}

\newcommand{\CN}[2]{\mathbb{C}\mathcal{N}\left(#1,#2\right)}
\newcommand{\vCN}[3]{\mathbb{C}\mathcal{N}_{#1}\left(#2,#3\right)}
\newcommand{\mCN}[4]{\mathbb{C}\mathcal{N}_{#1}\left(#2,#3,#4\right)}
\newcommand{\CW}[3]{\mathbb{C}\mathcal{W}_{#1}\left(#2,#3\right)}

\newcommand{\Cchisquare}[2]{\mathbb{C}\chi^{2}_{#1}(#2)}
\newcommand{\CF}[2]{\mathbb{C}\mathcal{F}_{#1}\left(#2\right)}
\newcommand{\betapdf}[1]{\mathcal{B}_{#1}}

\newcommand{\dist}{\overset{d}{=}}

\newcommand{\tAMF}{t_{\text{AMF}}}

\begin{document}
\title{Impact of covariance mismatched training samples on constant false alarm rate detectors}
\author{Olivier Besson \thanks{The author is with Universit\'{e} de Toulouse, ISAE-SUPAERO, 10 avenue Edouard Belin, 31055 Toulouse, France. Email: olivier.besson@isae-supaero.fr}}
\date{December 2020}
\maketitle
\begin{abstract}
The framework of this paper is that of adaptive detection in Gaussian noise with unknown covariance matrix when the training samples do not share the same covariance matrix as the vector under test. We consider a class of constant false alarm rate detectors which depend on two statistics $(\beta,\ttilde)$ whose distribution is parameter-free in the case of no mismatch and we analyze the impact of covariance mismatched training samples. More precisely, we provide a statistical representation of these two variables for an arbitrary mismatch. We show that covariance mismatch induces significant variations of the probability of false alarm and we investigate a way to  mitigate this effect.
\end{abstract}

\section{Introduction \label{sec:introduction}}
Deciding of the presence of a signal of interest (SoI) among Gaussian noise with unknown covariance matrix with the help of training samples that share the same covariance is a fundamental problem in multichannel processing, especially in radar applications where detection of a target with given space and/or time signature among clutter and noise is a key part of any radar architecture \cite{Richards10,Melvin12}. The papers and technical reports written by Kelly \cite{Kelly85,Kelly86,Kelly87,Kelly89,Kelly89b} have undoubtedly been the most influential work to tackle, to solve and to analyze this problem. His development and thorough analysis of the generalized likelihood ratio test (GLRT) through Wishart matrices theory have set the tempo and paved the way to numerous works around this fundamental theme, see e.g., \cite{Robey92,Kalson92,Bose95,Bose96b,Bose96,Raghavan95,Raghavan96,Kraut01,Pulsone00,Pulsone01,Richmond00,Conte01,Conte03,Jin05,Abramovich07,Bandiera09b,Bandiera08,Bandiera09,DeMaio16,Besson17} for a sample. When the SoI signature is known and the training samples have the same covariance matrix as the noise in the test vector, Bose and Steinhardt \cite{Bose95,Bose96b} showed that the maximal invariant statistic (MIS) for this detection problem is bi-dimensional and thus any detector which is a function of this 2-D maximal invariant statistic enjoys the constant false alarm rate (CFAR) property. The recent paper \cite{Coluccia20} provides a very clear review of CFAR detectors, introduces the so-called CFAR feature plane constructed from the MIS and investigates the use of  linear or non-linear detectors in this plane.

Deviations from this ideal situation often occur in practice. The first type of mismatch that has been dealt with is that on the SoI signature $\vv$, see e.g. \cite{Kelly89b,Bose96,Kalson92,Coluccia20}. Under this situation a major concern is to be able to maintain a constant false alarm rate while monitoring the trade-off between detection of slightly mismatched signals and rejection of unwanted signals. Ensuring the CFAR property when a mismatch about the SoI signature occurs is generally not a very complicated issue as the mismatch does not impact the distribution of the data and hence of the test statistic under the null hypothesis $H_{0}$. Therefore any detector which is CFAR when $\vv$ is known perfectly remains CFAR when a mismatch is present.

In contrast, a mismatch between the covariance matrix $\mSigmat$ of the training samples and the covariance matrix $\mSigma$ of the vector under test affects the distribution of any test statistic under $H_{0}$. Covariance mismatch can have a significant impact on most well-known detectors. In particular, it may preclude from maintaining a probability of false alarm $\Pfa$  that is independent of $\mSigma$ (and $\mSigmat$). In other words, a threshold which has been set to obtain a desired $\Pfabar$ when $\mSigma=\mSigmat$ results in an actual $\Pfa \neq \Pfabar$ when $\mSigma \neq \mSigmat$.  Unfortunately covariance mismatches can exist, starting with the  case of a partially homogeneous environment for which $\mSigmat = \gamma \mSigma$.  In this case the usual remedy is to apply some normalization so as to be insensitive to scaling, as is done by the adaptive coherence estimator \cite{Kraut99,Kraut05}. Some properties of the noise, e.g., that its power spectrum is symmetric can also be taken into account in order to improve performance \cite{Foglia17}. The reference \cite{Liu20} considers an extension to the case where noise power fluctuations between the range cells is observed and detection architectures are proposed to handle this situation. The training samples can also be contaminated by signal-like components or outliers \cite{Gerlach95,Gerlach02} or there might exist a  rank-one difference between $\mSigmat$ and $\mSigma$, for instance a surprise or undernulled interference  \cite{Richmond00c,Besson07c,Besson07g}.  

References \cite{Richmond00b,Blum00,McDonald00,Raghavan19} address rather general cases where the mismatch between $\mSigma$ and $\mSigmat$ is practically arbitrary, and they are closely related to the present work. In \cite{Richmond00b}, Richmond provides a thorough analysis of a large class of adaptive detectors in the case where $\mSigma$ and $\mSigmat$ satisfy the so-called generalized eigenrelation (GER), viz $\mSigmat^{-1}\vv = \lambda \mSigma^{-1}\vv$. Blum and McDonald \cite{Blum00} consider the case where a mismatch about the signal (target) signature is also present and a GER-type relation is satisfied, i.e. $\mSigmat^{-1}\vvbar = \lambda \mSigma^{-1}\vvbar$ where $\vvbar$ is the assumed value of $\vv$ used in the test statistic. They analyze Kalson's type of detectors \cite{Kalson92}. Probabilities of false alarm and probabilities of detection are obtained in integral form for any couple $(\mSigma,\mSigmat)$ satisfying the GER. The latter constraint is relaxed in \cite{McDonald00} where $\mSigma$ and $\mSigmat$ are arbitrary and where integral forms of $\Pfa$ and $\Pd$ are provided.  Recently Raghavan \cite{Raghavan19}  analyzed the performance of the adaptive matched filter (AMF) for arbitrary $\mSigma$ and $\mSigmat$.  The distribution of $\tAMF$ condition to the SNR loss is derived, and subsequently the corresponding false alarm rate and probability of detection.

In this paper, we consider an arbitrary mismatch between $\mSigma$ and $\mSigmat$ and we provide stochastic representations of two variables which form a maximal invariant statistic in the absence of mismatch. The new stochastic representations are given in terms of well-known distributions and enable one to quickly grasp how the maximal invariant statistic is impacted. They generalize the expressions obtained in \cite{Richmond00b} where the GER was enforced. In contrast to \cite{Blum00,McDonald00} we do not address the problem of deriving an expression for $\Pfa$ and $\Pd$ here, which anyhow can be obtained in integral form only. The objective of this paper is more  to study how $\Pfa$ deviates from $\Pfabar$ (calculated when $\mSigmat=\mSigma$)  when $\mSigmat$ varies randomly around $\mSigma$ and to investigate possible ways to mitigate the fluctuations caused by covariance mismatch. 

The paper is organized as follows. In Section \ref{sec:analysis} we derive a stochastic representation of random variables $\beta$ and $\ttilde$ (to be defined below) for arbitrary $\mSigma$ and $\mSigmat$. This representation has a very similar form as in the no mismatch case and provides quick insights into  how the distribution of these variables is affected. When the generalized eigenrelation is assumed to hold our representation boils down to that of \cite{Richmond00b}. Moreover we show that under the GER a clairvoyant CFAR detector is possible, which depends on a single scalar parameter namely $(\vv^{H}\mSigmat^{-1}\vv)/(\vv^{H}\mSigma^{-1}\vv)$. Based on the form of this clairvoyant detector, we propose in Section \ref{sec:mitigation} a solution to somehow mitigate the effects of covariance mismatch and to render the false alarm rate less sensitive to $\mSigmat \neq \mSigma$. Finally we investigate how $\Pd$ is impacted.

\emph{Notations}: Before proceeding, we state some notations regarding the distributions used frequently in the paper. The notation $A \dist B$ means that the random variables $A$ and $B$ have the same distribution. The  $N$-dimensional complex multivariate Gaussian distribution is denoted as $\vx \dist \vCN{N}{\vxbar}{\mSigma}$ and its probability density function (p.d.f.) is given by $p(\vx) = \pi^{-N} \det{\mSigma}^{-1} \etr{-(\vx-\vxbar)^{H} \mSigma^{-1} (\vx-\vxbar)}$.  The extension to a complex matrix-variate Gaussian distribution for a $N \times K$ matrix $\X$ is denoted as $\X \dist \mCN{N,K}{\Xbar}{\mSigma}{\mPsi}$ with density $p(\X) = \pi^{-NK} \det{\mSigma}^{-K} \det{\mPsi}^{-N}\etr{-(\X-\Xbar)^{H}\mSigma^{-1}(\X-\Xbar)\mPsi^{-1}}$. If $\vx \dist \vCN{p}{\vxbar}{\eye{p}}$  then $\vx^{H}\vx$ follows a complex chi-square distribution with $p$ degrees of freedom and non centrality parameter  $\delta =  \vxbar^{H}\vxbar$, i.e.  $\vx^{H}\vx \dist \Cchisquare{p}{\delta}$.  If $\X \dist \mCN{N,K}{\Mzero}{\mSigma}{\eye{K}}$ then $\W = \X \X^{H} \dist \CW{N}{K}{\mSigma}$ follows a complex Wishart distribution with parameters $K$ and $\mSigma$.   Let $U \dist \Cchisquare{p}{\delta}$ and $V \dist \Cchisquare{q}{0}$. Then $F = U/V$ follows a complex F-distribution, which is noted as  $F \dist \CF{p,q}{\delta}$.  When $\delta=0$,   $B=(1+F)^{-1} \dist \betapdf{p,q}$ has a central beta distribution.

For any rectangular matrix $\Z$, whenever it is partitioned as $\Z = \begin{pmatrix} \Z_{1} \\ \Z_{2} \end{pmatrix}$ then $ \Z_{1}$ contains the first $(N-1)$ rows and all columns of $\Z$. If a square matrix $\mOmega$ is partitioned as $\mOmega = \begin{pmatrix} \mOmega_{11}  & \mOmega_{12} \\ \mOmega_{21} & \Omega_{22} \end{pmatrix}$ then $\mOmega_{11} $ is a $(N-1)\times (N-1)$ matrix.

\section{Analysis of the impact of $\mSigmat \neq \mSigma$ \label{sec:analysis}}
As stated above, we consider the generic composite hypothesis testing problem
\begin{align}\label{prob_detect}
	H_{0}: \, &\vx \dist \vCN{N}{\vzero}{\mSigma}; \Xt \dist \mCN{N,K}{\Mzero}{\mSigmat}{\eye{K}} \nonumber \\
	H_{1}: \, &\vx \dist \vCN{N}{\alpha \vv}{\mSigma}; \Xt \dist \mCN{N,K}{\Mzero}{\mSigmat}{\eye{K}}
\end{align}
We begin by recalling some important results in the case where $\mSigmat=\mSigma$ then address the case $\mSigmat \neq \mSigma$. 
\subsection{Background on the case of no mismatch}
When $\mSigmat=\mSigma$  \cite{Bose95,Bose96b} showed that the maximal invariant statistic for the detection problem in \eqref{prob_detect} is 2-D and is a bijective function of $(s_{1},s_{2})$ defined as
\begin{equation}
s_{1}  =\vx^{H}\St^{-1}\vx ; \, s_{2} = \frac{|\vx^{H}\St^{-1}\vv|^{2}}{\vv^{H}\St^{-1}\vv}
\end{equation}
where $\St=\Xt\Xt^{H}$ denotes the sample covariance matrix (SCM) of the training samples. Bose and Steinhardt used as the maximal invariant statistic $\rho=(1+s_{1}-s_{2})^{-1}$ and $\eta=(1+s_{1})^{-1}(1+s_{1}-s_{2})$ where $\rho$ is the loss factor and $\eta$ is related to Kelly's test statistic. The joint distribution of $(\rho,\eta)$ was derived under both hypotheses $H_{0}$ and $H_{1}$. By integrating over this distribution \cite{Bose96b} manages to obtain the probability of false alarm ($\Pfa$) and probability of detection ($\Pd$) of the locally uniformly most powerful invariant test. Note that any function of  $(\rho,\eta)$ yields a CFAR detector. In the sequel similarly to \cite{Richmond00,Coluccia20} we choose to make use  of the pair of random variables $(\beta,\ttilde)$ (also a MIS) defined as
\begin{equation}
\beta = \frac{1}{1+s_{1}-s_{2}}	; \,\ttilde = \frac{s_{2}}{1+s_{1}-s_{2}}
\end{equation}
$\beta$ corresponds to the loss factor while $\ttilde$  is the test statistic for Kelly's GLRT.  It has been shown, see e.g.,  \cite{Kelly86,Bose96b,Richmond00}, that when $\mSigma=\mSigmat$,
\begin{align}
\beta & \dist \left[1 + \frac{\Cchisquare{N-1}{0}}{\Cchisquare{K-N+2}{0}}\right]^{-1} \dist \betapdf{N-1,K-N+2} \label{pdf_beta_no_mismatch} \\
\ttilde | \beta & \dist \CF{1,K-N+1}{\beta |\alpha|^{2}\vv^{H}\mSigma^{-1}\vv} \label{pdf_ttilde_no_mismatch}
\end{align}
As a consequence, under $H_{0}$, $\beta$ and $\ttilde$ are independent and $\ttilde \dist \CF{1,K-N+1}{0}$. Therefore, the threshold of any detector which is a function of $(\beta$, $\ttilde)$ can be set to ensure a given probability of false alarm whatever $\mSigma$.

\subsection{Representation of $\beta$ and $\ttilde$ in case of covariance mismatch}
We now aim at obtaining a representation of $(\beta,\ttilde)$ equivalent to \eqref{pdf_beta_no_mismatch}-\eqref{pdf_ttilde_no_mismatch} for $\mSigmat \neq \mSigma$. Prior to that we mention that the recent paper \cite{Raghavan19} analyzes the probability of false alarm and probability of detection of the AMF in case of arbitrary mismatch between $\mSigma$ and $\mSigmat$. In \cite{Raghavan19}  the conditional distribution of $\tAMF$ given the SNR loss $\rho$ is derived and conditional and unconditional $\Pfa$ and $\Pd$ are obtained. The derivations used below to obtain a representation of $(\beta,\ttilde)$ borrow from techniques used in \cite{Raghavan19} and many papers, viz. a whitening step, followed by rotation and finally the use of partitioned Wishart matrices properties. The main result is stated below.

When $\mSigmat \neq \mSigma$ $\beta$ and $\ttilde$ have the following stochastic representation:
\begin{align}
&\beta  \dist (1+\vxtilde_{1}\W_{11}^{-1}\vxtilde_{1})^{-1} \dist \left[1+\frac{\sum_{i=1}^{N-1}\lambda_{i} \Cchisquare{1}{0}}{\Cchisquare{K-N+2}{0}}\right]^{-1}  \label{storep_beta_mismatch} \\
&\ttilde | \vxtilde_{1},\W_{11}  \dist \left[1 + \beta \left(\frac{\vv^{H}\mSigmat^{-1}\vv}{\vv^{H}\mSigma^{-1}\vv}-1\right)\right] \nonumber \\
&\, \times \CF{1,K-N+1}{\beta \frac{ ||\alpha| (\vv^{H}\mSigmat^{-1}\vv)^{1/2}+\mOmega_{21}\mOmega_{11}^{-1}\vxtilde_{1}|^{2}}{1+ \beta \left( \frac{\vv^{H}\mSigmat^{-1}\vv}{\vv^{H}\mSigma^{-1}\vv}-1\right)}} \label{storep_ttilde_mismatch}
\end{align}

where $\W_{11} \dist \CW{N-1}{K}{\eye{N-1}}$ is independent of  $\vxtilde_{1} \dist \vCN{N-1}{\vzero}{\mOmega_{11}}$. In the previous equation, 
\begin{align}
\mOmega &= \Q^{H}\mSigmat^{-1/2} \mSigma \mSigmat^{-H/2} \Q \nonumber \\
&= \begin{pmatrix} \underset{N-1|N-1}{\mOmega_{11}}  & \underset{N-1|1}{\mOmega_{12}} \\ \underset{1|N-1}{\mOmega_{21}} & \underset{1|1}{\Omega_{22}} \end{pmatrix}
\end{align}
and  $\lambda_{i}$  are the eigenvalues of $\mOmega_{11}$. $\Q$ is a unitary matrix which rotates the whitened vector $\mSigmat^{-1/2}\vv$, i.e. $\Q^{H}\mSigmat^{-1/2}\vv = (\vv^{H}\mSigmat^{-1}\vv)^{1/2}\elast$ with $\elast = \begin{bmatrix} \vzero^{T} & 1 \end{bmatrix}^{T}$.

The proof of this result is given in Appendix \ref{app:proof}. \textit{Equations \eqref{storep_beta_mismatch}-\eqref{storep_ttilde_mismatch} provide the statistical representation of $\beta$ and $\ttilde$ for any couple of matrices $(\mSigma,\mSigmat)$ and whatever the value of $\alpha$, i.e., under $H_{0}$ or $H_{1}$.} 

Some comments are in order. First of all, \emph{in the case of no mismatch} $\mOmega=\eye{N}$, $\lambda_{i}=1$, $\mOmega_{21}=\mat{0}$, and \emph{equations \eqref{storep_beta_mismatch}-\eqref{storep_ttilde_mismatch} reduce to equations \eqref{pdf_beta_no_mismatch}-\eqref{pdf_ttilde_no_mismatch}}. 

Next, while there is quite a striking similarity between \eqref{pdf_beta_no_mismatch}-\eqref{pdf_ttilde_no_mismatch} and \eqref{storep_beta_mismatch}-\eqref{storep_ttilde_mismatch},  one can observe \textbf{three big differences} in the case $\mSigmat \neq \mSigma$:
\begin{enumerate}
\item the random variable \emph{$\beta$ is no longer beta distributed} since the term $\Cchisquare{N-1}{0}$ is replaced by $\sum_{i=1}^{N-1}\lambda_{i} \Cchisquare{1}{0}$ which corresponds to a quadratic form in i.i.d. central complex normal variables. As shown in \cite{Besson20e}, the distribution of $\beta$ can be well approximated by $a\Cchisquare{\nu}{0}+b$ or even by $a'\Cchisquare{\nu'}{0}$. 
\item as for $\ttilde$, \emph{the coefficient before the F distribution is no longer $1$} but
\begin{equation*}
1 + \beta \left(\frac{\vv^{H}\mSigmat^{-1}\vv}{\vv^{H}\mSigma^{-1}\vv}-1\right)
\end{equation*} 
Therefore, depending on $(\vv^{H}\mSigmat^{-1}\vv)/(\vv^{H}\mSigma^{-1}\vv)$, this coefficient can be above or under unity.
\item \emph{the complex F distribution is no longer central under $H_{0}$}. In particular, it means that this random variable is likely to take on larger values than in the case of no mismatch, which implies an increase in the false alarm rate if $\ttilde$ is used as a test statistic, which is the case for Kelly's detector.
\end{enumerate}

Another important comment when examining \eqref{storep_ttilde_mismatch} is the importance of the term $\mOmega_{21}\mOmega_{11}^{-1}\vxtilde_{1}$. Indeed, if this term is zero, then the complex F distribution is central under $H_{0}$, as in the case of no mismatch. This happens if $\mOmega_{21}=\vzero$ which is equivalent to
\begin{align}\label{GER}
\mOmega_{21} = \vzero &\Leftrightarrow \vv^{H} \mSigmat^{-1}\mSigma \Vorth = \vzero \nonumber \\
&\Leftrightarrow  \mSigma\mSigmat^{-1}\vv = \lambda \vv \nonumber \\
&\Leftrightarrow  \mSigmat^{-1}\vv = \lambda \mSigma^{-1}\vv 
\end{align}
The previous equation corresponds to the so-called \textit{ generalized eigenrelation} \cite{Richmond00c}. If the GER is satisfied, then we have the simpler representation
\begin{align}\label{storep_ttilde_mismatch_GER}
\ttilde | \beta &\dist \left[1 + \beta \left(\frac{\vv^{H}\mSigmat^{-1}\vv}{\vv^{H}\mSigma^{-1}\vv}-1\right)\right] \nonumber \\
&\times\CF{1,K-N+1}{\beta \frac{ |\alpha|^{2} (\vv^{H}\mSigmat^{-1}\vv)}{1+ \beta \left( \frac{\vv^{H}\mSigmat^{-1}\vv}{\vv^{H}\mSigma^{-1}\vv}-1\right)}} 
\end{align}
which coincides with the expression in \cite{Richmond00b}. Under $H_{0}$ \eqref{storep_ttilde_mismatch_GER} reduces to
\begin{equation}
\ttilde | \beta ; H_{0} \dist \left[1 + \beta \left(\frac{\vv^{H}\mSigmat^{-1}\vv}{\vv^{H}\mSigma^{-1}\vv}-1\right)\right] \CF{1,K-N+1}{0}
\end{equation}
If $\vv^{H}\mSigmat^{-1}\vv=\vv^{H}\mSigma^{-1}\vv$ then $\ttilde$ is independent of $\beta$ and its distribution under $H_{0}$ does not longer depend on any parameter, which means that Kelly's detector would still be CFAR despite the covariance mismatch. Actually, this is to be related with the result in \cite{Besson07g} where one considers a surprise interference, i.e., $\mSigma=\mSigmat + \vq \vq^{H}$ where $\vq$ satisfies $\vq^{H}\mSigma^{-1}\vv=0$ for the GER to be true. In this case we have actually $\vv^{H}\mSigmat^{-1}\vv=\vv^{H}\mSigma^{-1}\vv$  and it was shown in \cite{Besson07g} that the GLRT is simply Kelly's detector. The analysis in this paper proves that Kelly's detector is CFAR in this situation.

Before closing this section, we note that \eqref{storep_ttilde_mismatch} and \eqref{storep_ttilde_mismatch_GER} provide the conditional distribution of $\ttilde$. Obtaining the probability of false alarm and probability of detection of Kelly's detector thus requires to integrate this conditional distribution. As has already been observed in \cite{Richmond00b,Blum00,McDonald00,Raghavan19} integration does not unfortunately lead to a closed-form expression. Therefore in the sequel the stochastic representations of $\beta$ in \eqref{storep_beta_mismatch} and $\ttilde$ in \eqref{storep_ttilde_mismatch} or \eqref{storep_ttilde_mismatch_GER} will be used to generate random samples and estimate $\Pfa$ or $\Pd$ from Monte-Carlo simulations.

\begin{figure*}[htb]
	\centering
	\subfigure{%
		\includegraphics[width=7.5cm]{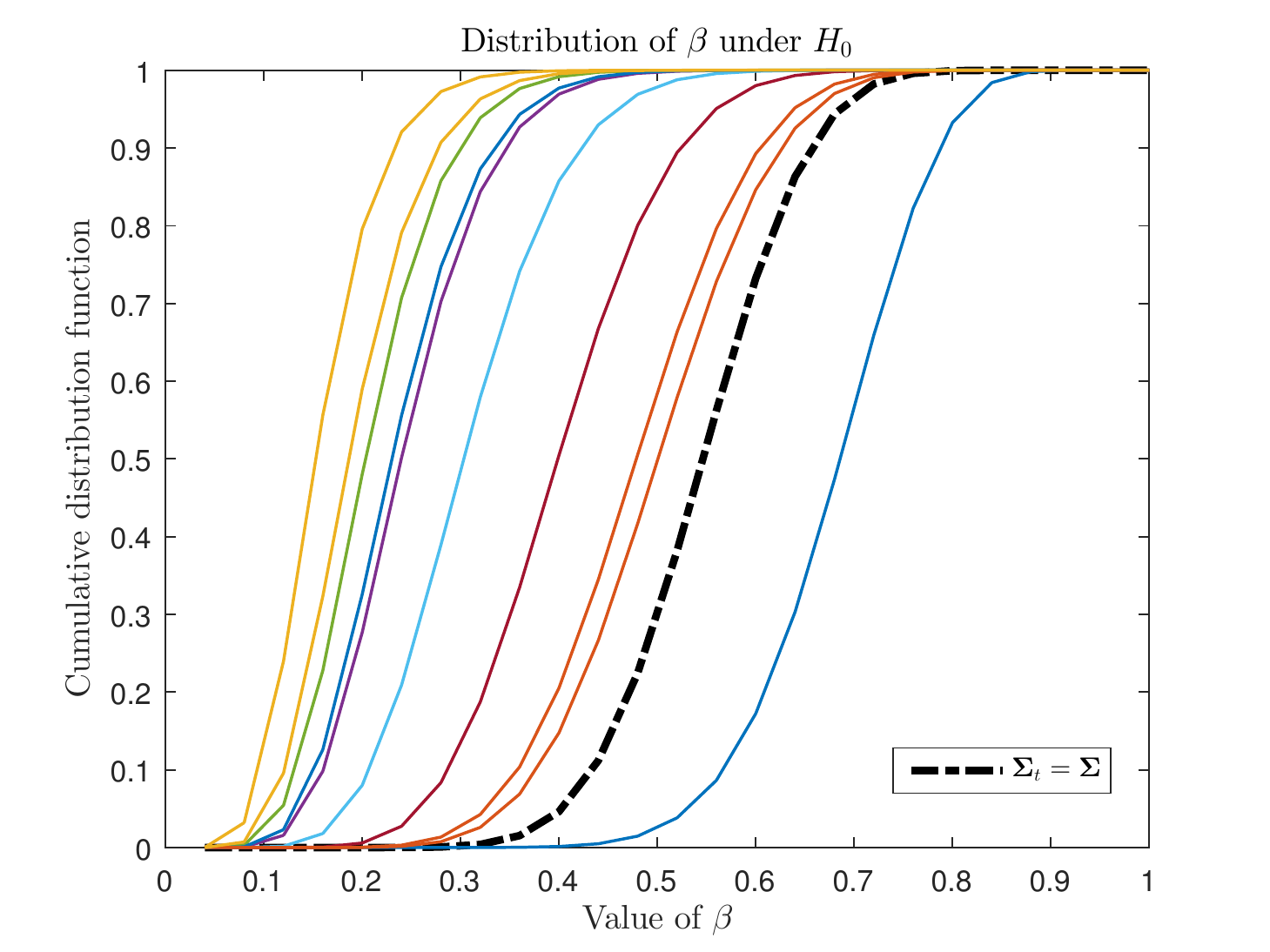}}	
	\subfigure{%
		\includegraphics[width=7.5cm]{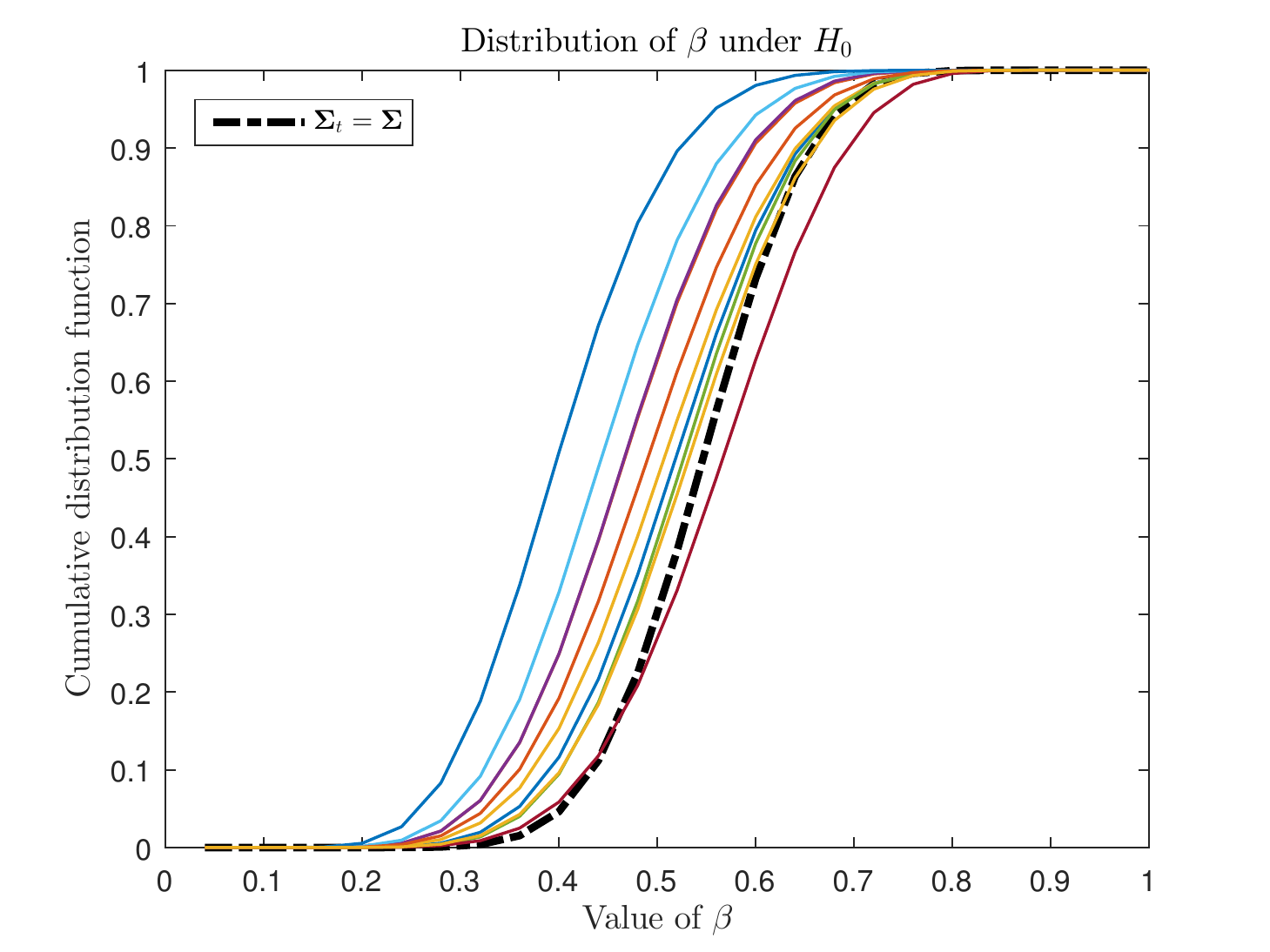}}\\ 	
	\subfigure{
		\includegraphics[width=7.5cm]{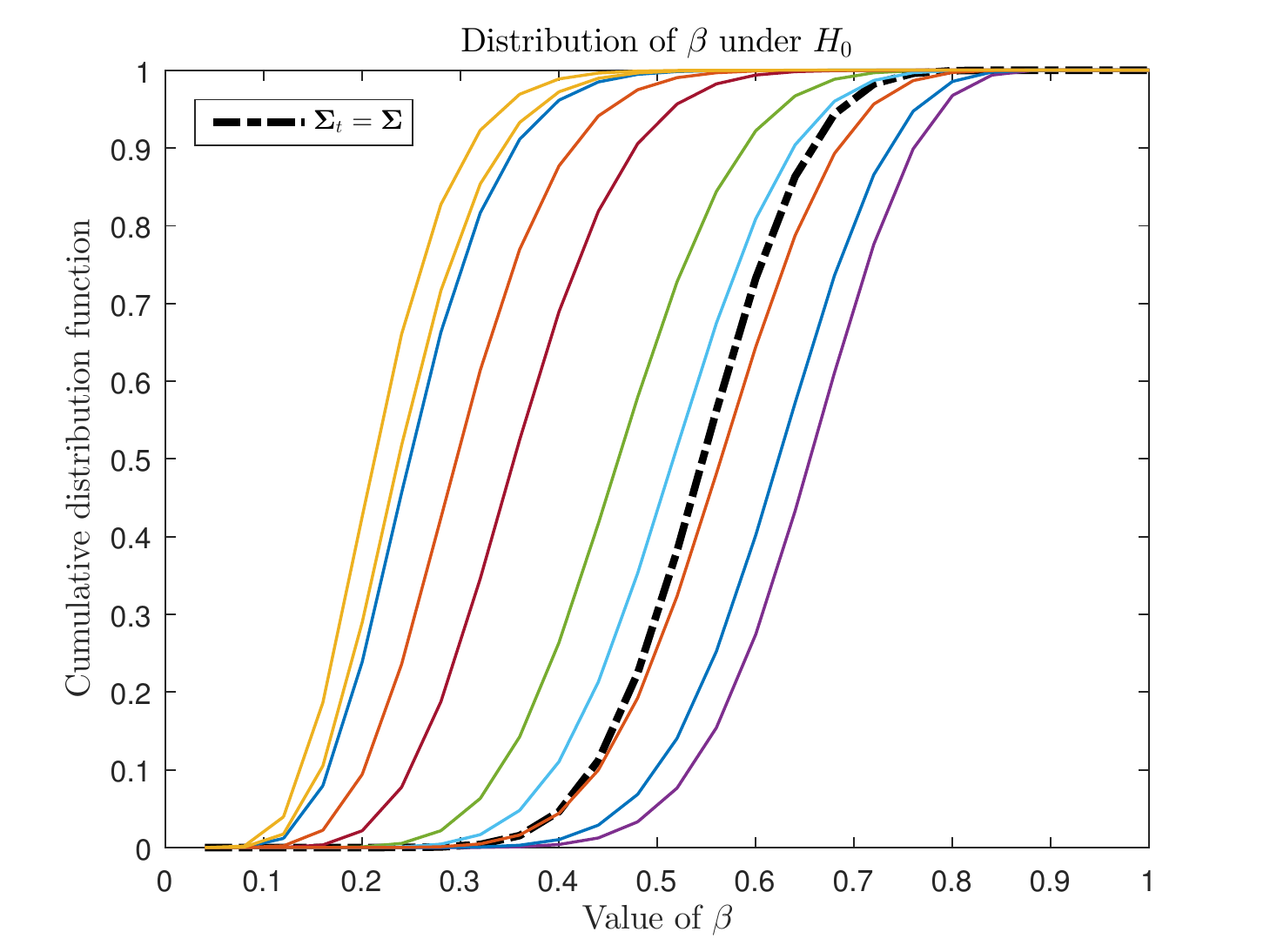}}
	\subfigure{%
		\includegraphics[width=7.5cm]{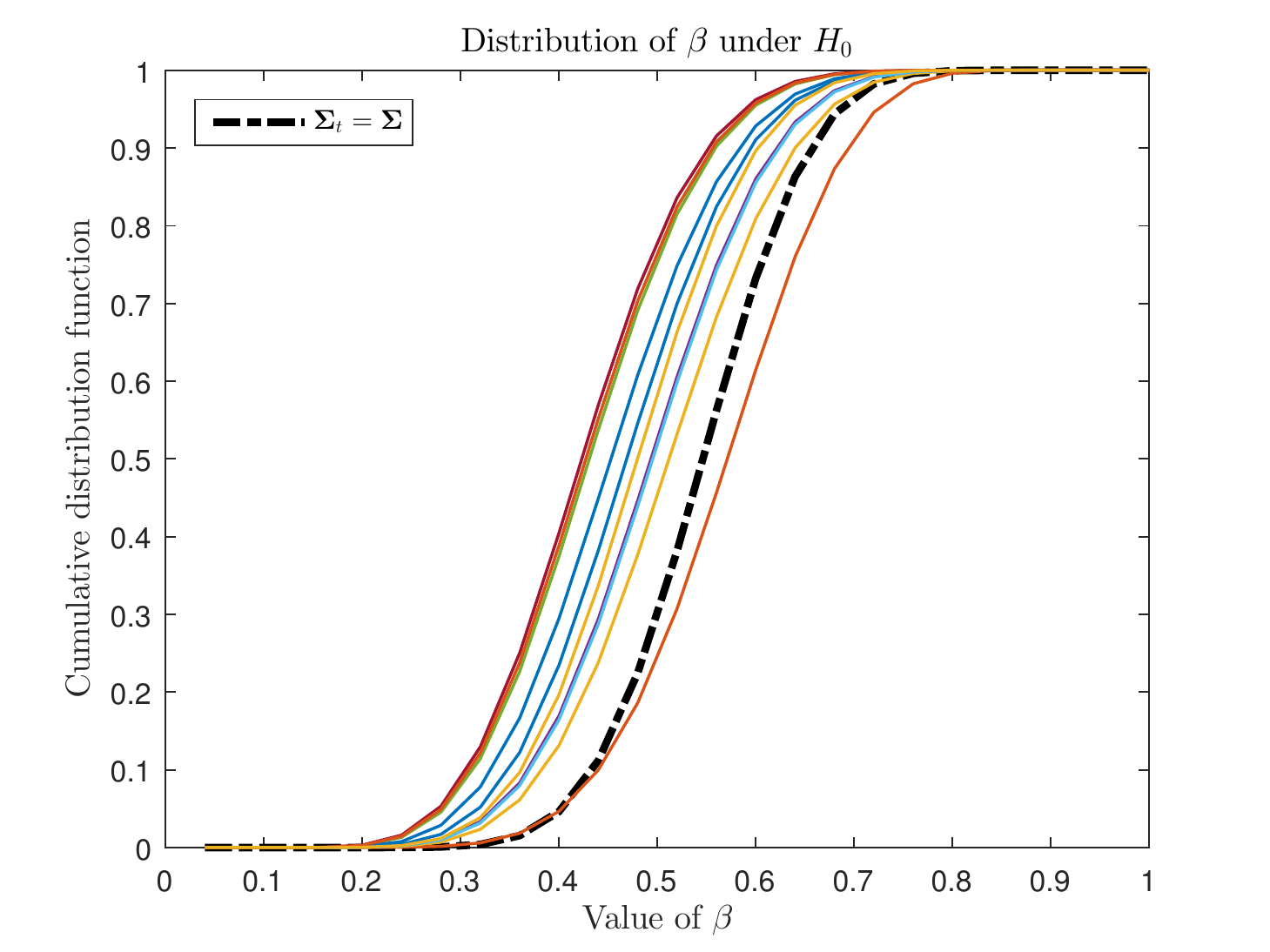}}
	\caption{Distribution of $\beta$  in the case of covariance mismatch  under $H_{0}$. The left panel concerns case 1, the right panel case 2. The lower panel corresponds to the case where the GER is satisfied.}
	\label{fig:cdf_beta_K=32_H0}
\end{figure*}

\begin{figure*}[htb]
	\centering
	\subfigure{%
		\includegraphics[width=7.5cm]{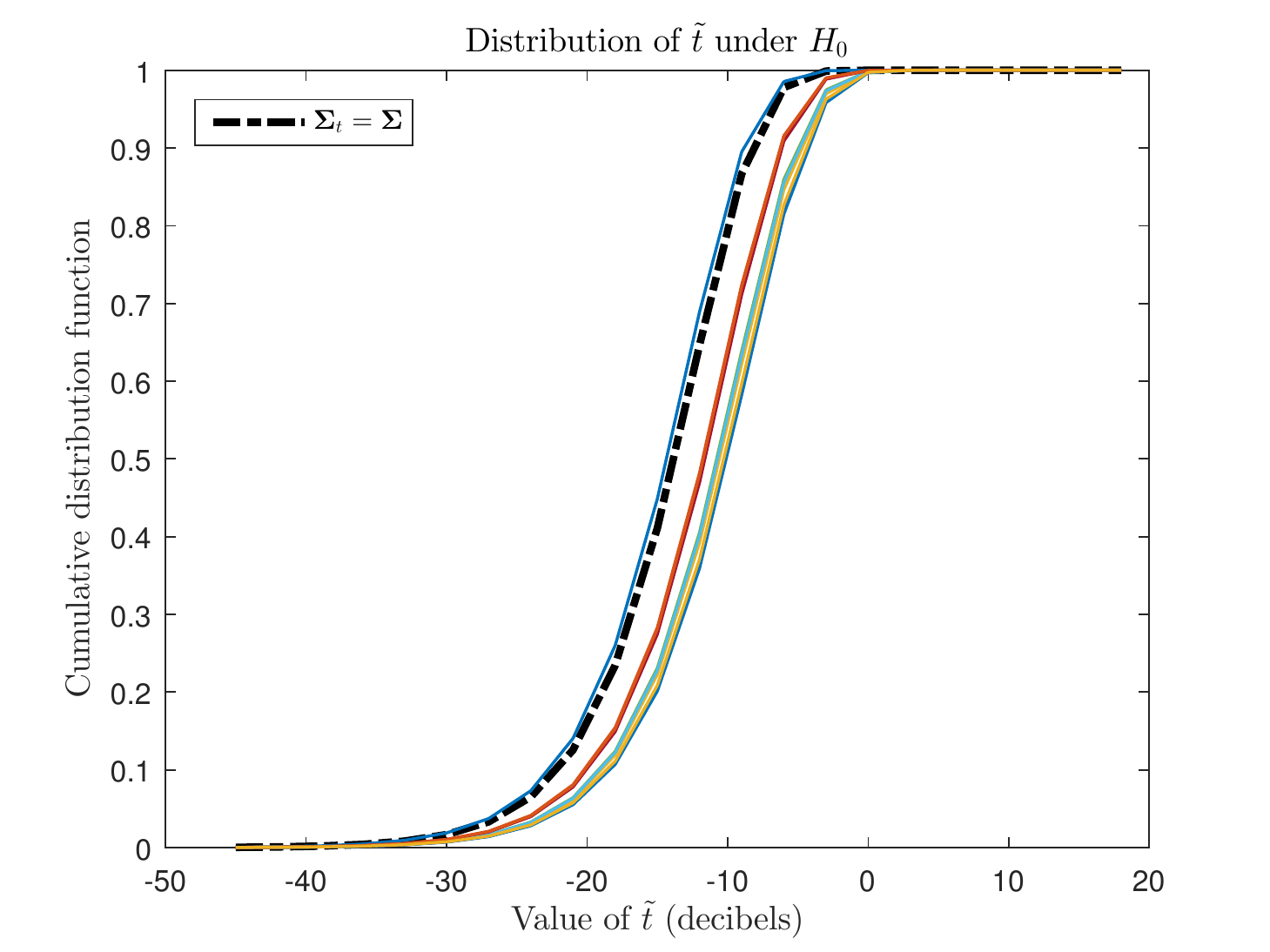}}
	\subfigure{%
		\includegraphics[width=7.5cm]{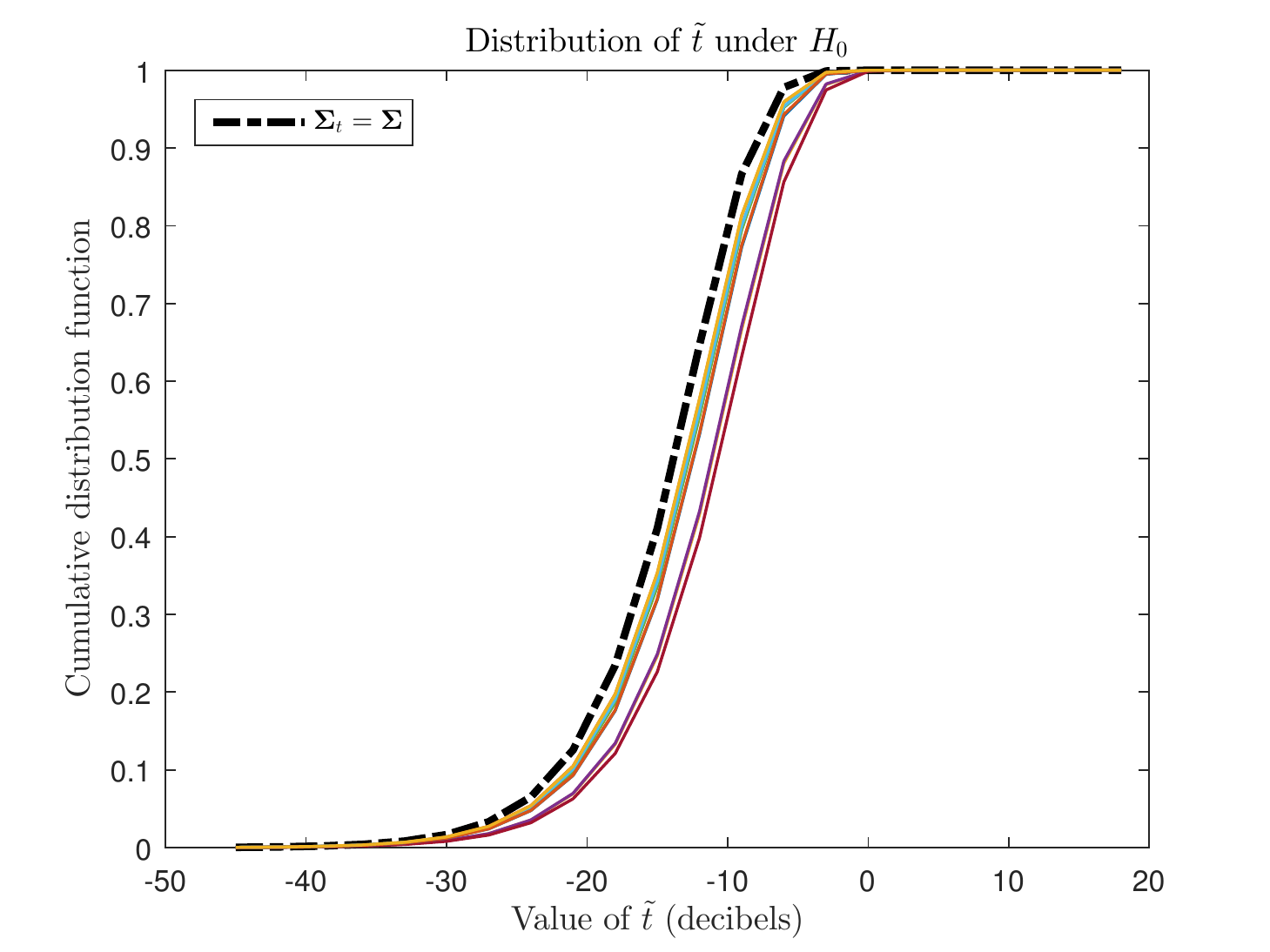}}\\ 		
	\subfigure{
		\includegraphics[width=7.5cm]{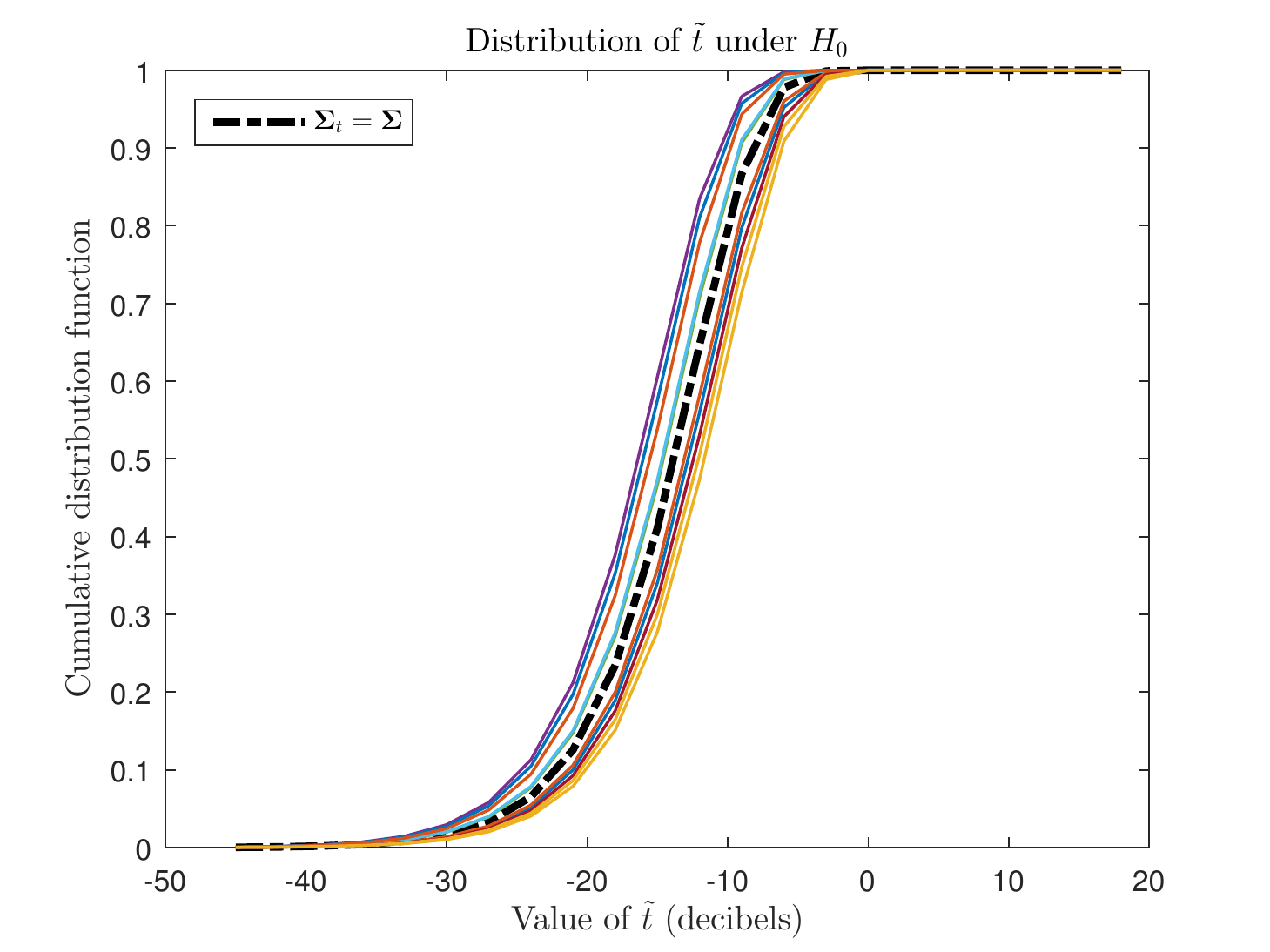}}
	\subfigure{%
		\includegraphics[width=7.5cm]{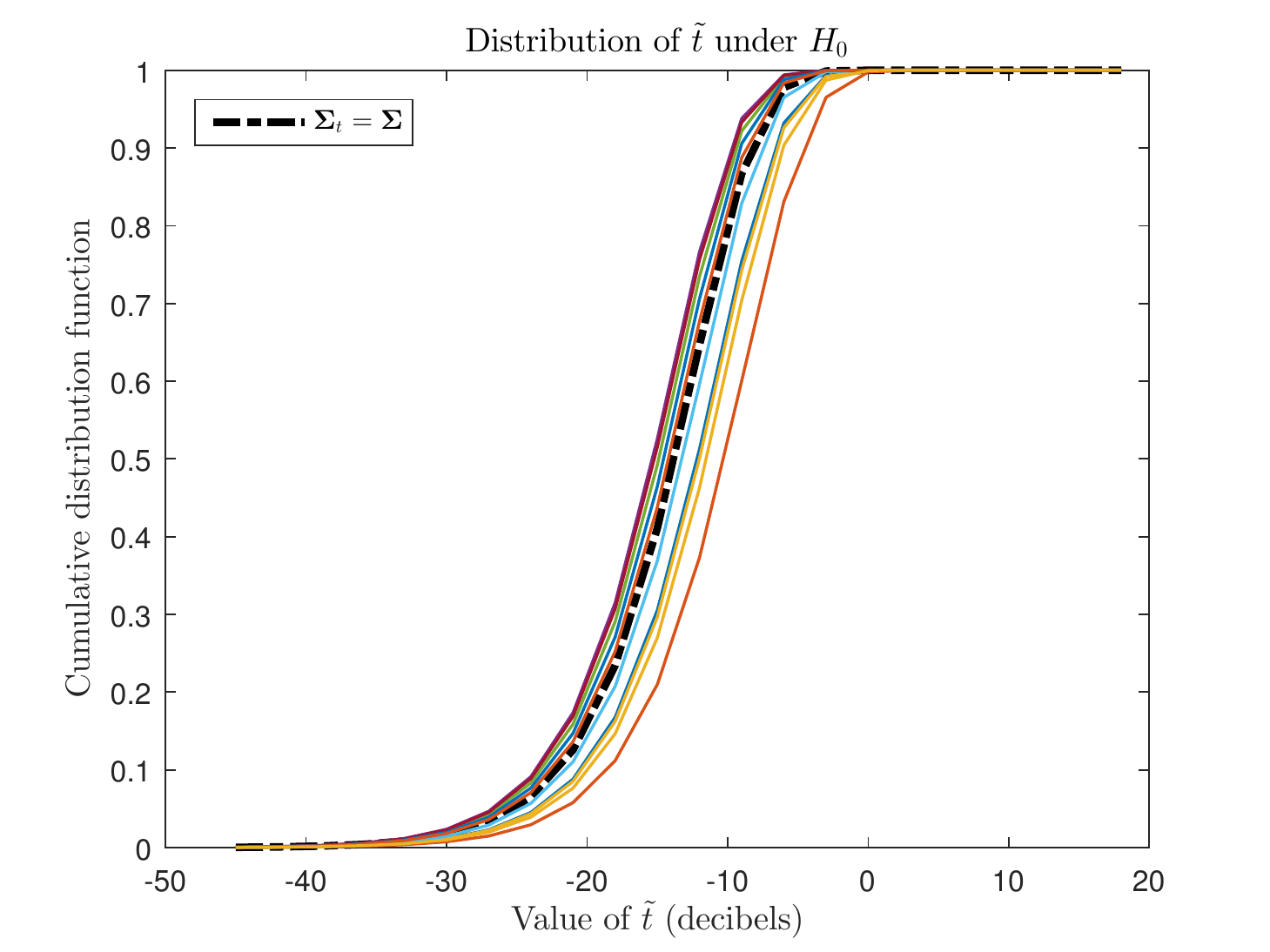}}
	\caption{Distribution of $\ttilde$ in the case of covariance mismatch under $H_{0}$. The left panel concerns case 1, the right panel case 2. The lower panel corresponds to the case where the GER is satisfied.}
	\label{fig:cdf_F_K=32_H0}
\end{figure*}

\subsection{Illustrations}
In this section, we provide illustrations of the theoretical results presented above. We consider a scenario with $N=16$. The noise covariance matrix is chosen as $\mSigma = P_{c} \mSigma_{c} + \eye{N}$ where the first component corresponds to clutter while the second is thermal white noise. We have $\mSigma_{c}(m,n) = \exp \{-2\pi^{2}\sigma_{f}^{2}(m-n)^{2}\}$ and $\sigma_{f}$ is such that the one-lag correlation is equal to $0.95$. The SoI signature is $\vv = N^{-1/2} \begin{bmatrix} 1 & e^{i2\pi f_d} & \ldots & e^{i2\pi (N-1) f_d} \end{bmatrix}^{T}$ with $f_{d}=0.08$ so that the target is slowly moving and buried in clutter. The clutter to noise ratio is $20$dB. The number of training samples is set to $K=2N$.

We will consider two different covariance mismatches. In the first case, $\mSigmat = \mSigma^{1/2} \Wt^{-1} \mSigma^{H/2}$ where $\Wt$ is drawn from a $\CW{N}{\nu}{\mu^{-1}\eye{N}}$ distribution. The mean value of $\Wt^{-1}$ is $\E{\Wt^{-1}}=\gamma\eye{N}$ with $\gamma=(\nu-N)^{-1}\mu$.  This enables to have a mean value $\E{\mSigmat} = \gamma \mSigma$ and $10\log_{10} \gamma$ is drawn from a uniform distribution in [$-6$dB, $6$dB].  This means that $\mSigmat$ can be quite different from $\mSigma$.   In the second case, $\mSigma$ and $\mSigmat$ share the same eigenvectors but different eigenvalues, i.e., $\mSigma=\U\mLambda\U^{H}$,  $\mSigmat=\U\mLambda^{1/2}\diag{\gamma_{n}}\mLambda^{1/2}\U^{H}$ and $10\log_{10}\gamma_{n}$ is drawn uniformly in [$-6$dB, $6$dB]. For both cases, the GER does not hold. In order to obtain similar types of mismatch but with the GER satisfied we can proceed as follows. In case 1, as shown in \cite{Besson20e}, the GER is satisfied when $\mSigmat$ is  generated $\mSigmat$ as
\begin{equation}\label{GER_chol}
\Qv^{H} \mSigmat \Qv = \chol{\Qv^{H} \mSigma \Qv} \begin{pmatrix} \mPsi_{11}^{-1} & \Mzero \\ \vzero & \Psi_{22}^{-1} \end{pmatrix} \chol{\Qv^{H} \mSigma \Qv}
\end{equation}
where $\Qv=\begin{bmatrix} \Vorth & \vv \end{bmatrix}$ and $\chol{.}$ stands for the Cholesky factor. $\mPsi_{11}$ will be drawn from a Wishart distribution and $\Psi_{22}$ from a chi-square distribution with $\E{\mPsi_{11}^{-1}}=\gamma \eye{N-1}$ and $\E{\Psi_{22}^{-1}}=\gamma$. Again,  $10\log_{10}\gamma$ is drawn uniformly in [$-6$dB, $6$dB]. As for case 2 with eigenvalues mismatch, we show in the appendix that $\mSigmat=\U\mLambda^{1/2}\Wt\mLambda^{1/2}\U^{H}$  with $\Wt$ given by \eqref{gen_GER_wLambda} enables one to satisfy the GER. As before, the eigenvalues $10\log_{10}\vl_{1}(n)$ and $10\log_{10}\ell_{2}$ are drawn uniformly in [$-6$dB, $6$dB]. We thus have two kinds of mismatches and, for each case, a variant that enforces the GER.

First, we investigate how the distributions of $\beta$ and $\ttilde$ are affected. Towards this end, we generated $10$ random matrices  $\mSigmat$ and we display the distributions  of $\beta$ and $\ttilde$ under $H_{0}$ for these $10$ different matrices $\mSigmat$. For comparison purposes we also display (dashed-dotted lines) their distributions in the case of no mismatch.  The results are shown in Figure \ref{fig:cdf_beta_K=32_H0} for $\beta$ and Figure \ref{fig:cdf_F_K=32_H0} for $\ttilde$. The major observations to be made are the following:
\begin{itemize}
	\item $\beta$ most frequently takes lower values when $\mSigmat \neq \mSigma$ than when $\mSigmat=\mSigma$, although the converse may be true in some instances. Then, for the AMF whose test statistics is $\ttilde / \beta$, it means that the AMF value will likely be larger than expected, causing an increase in $\Pfa$.
	\item we observe that $\beta$ undergoes more important variations in case 1 than in case 2 (mismatched eigenvalues) where the variations are less significant. Generally the mismatch of case 1, of the form $\mSigmat = \mSigma^{1/2} \Wt^{-1} \mSigma^{H/2}$, seems to be more severe.
	\item $\ttilde$ tends to increase when $\mSigmat \neq \mSigma$ compared to the case $\mSigmat=\mSigma$, especially when the GER is not satisfied. It seems that the difference is less pronounced when the GER is enforced. An increase in $\ttilde$ will result in an increase of the false alarm rate of Kelly's detector.
\end{itemize}

\begin{figure}[htb]
\centering
\subfigure{%
\includegraphics[width=8cm]{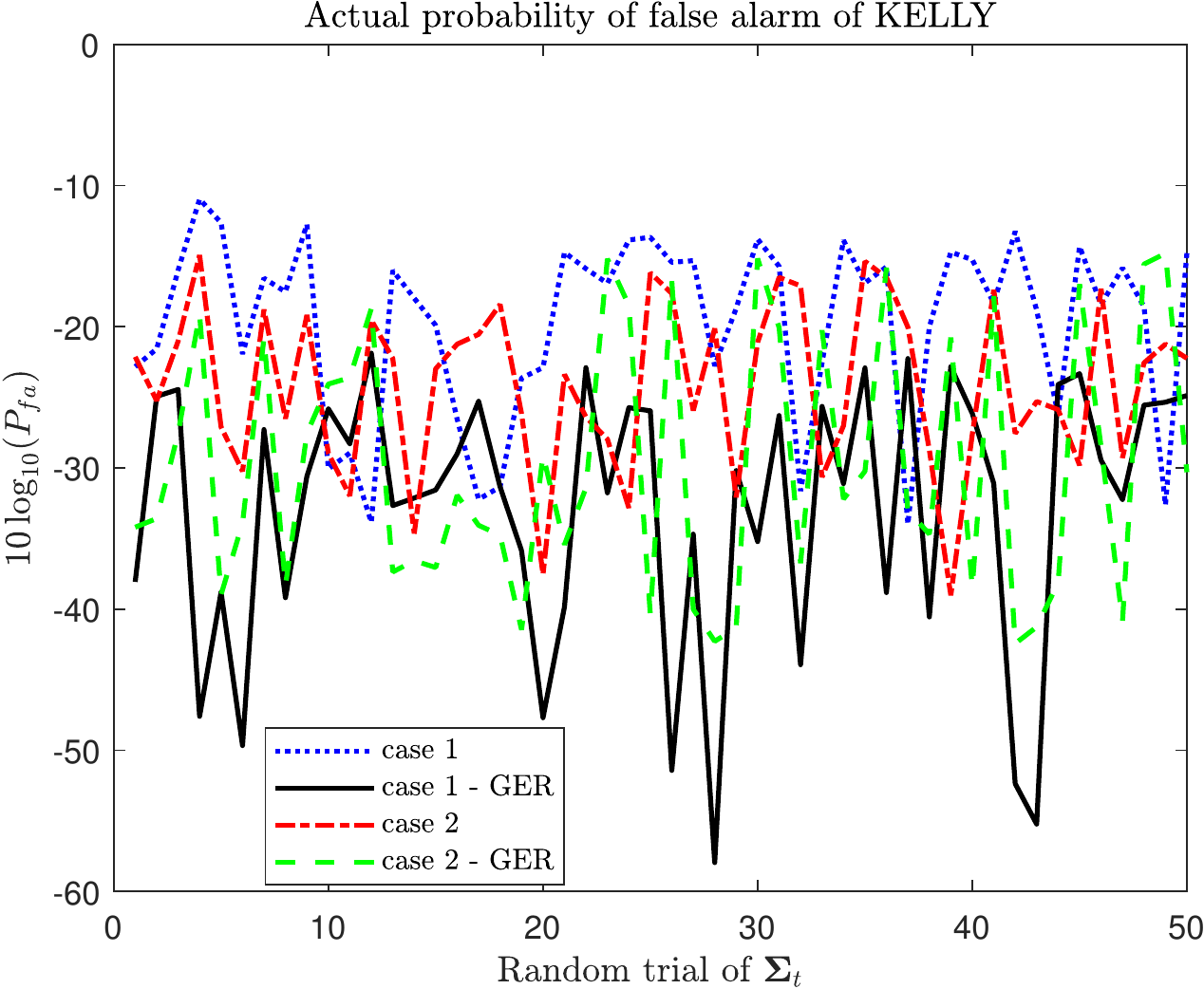}} \\
\subfigure{%
	\includegraphics[width=8cm]{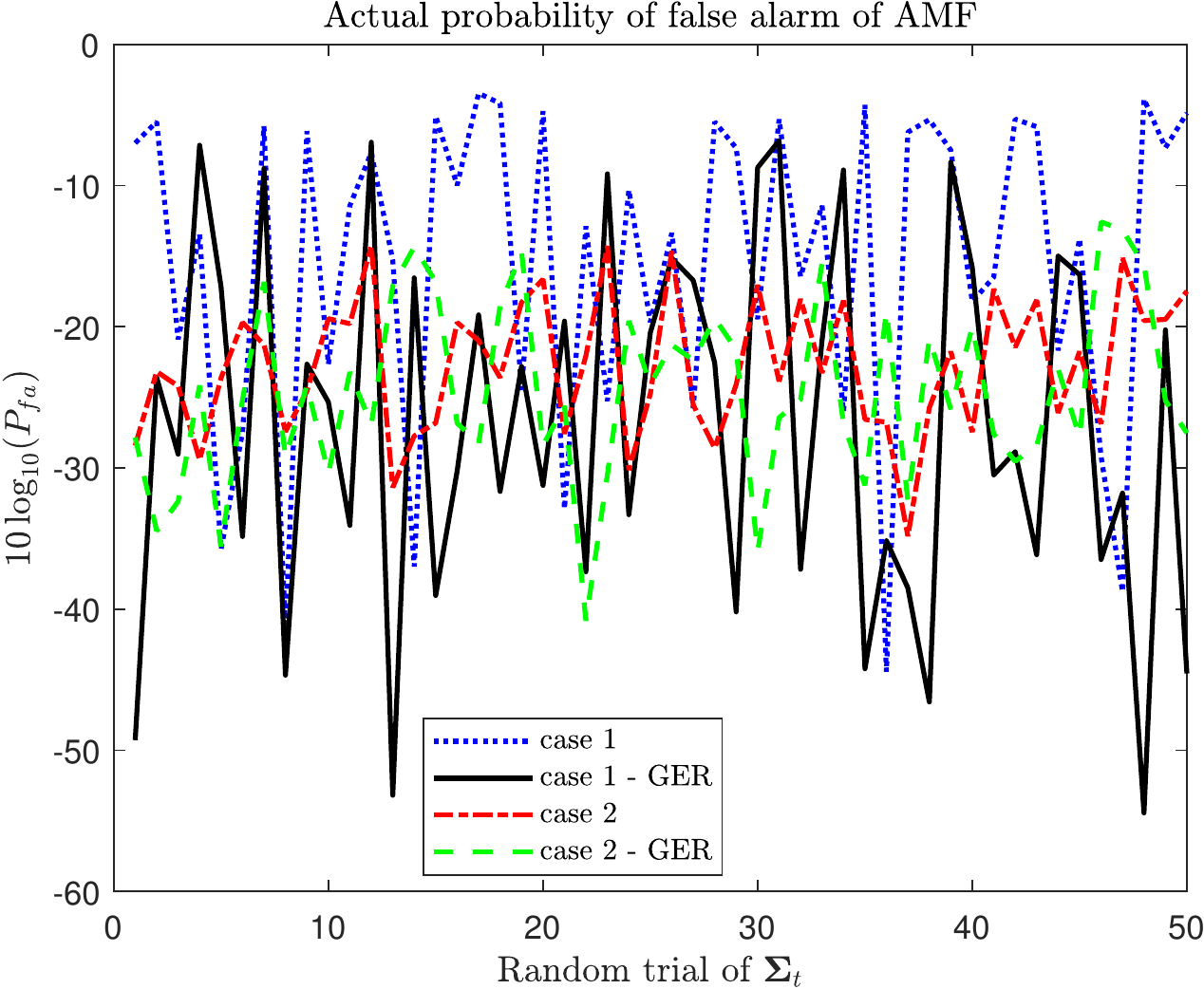}}
\caption{Actual probability of false alarm of Kelly ($\ttilde$) and AMF detector ($\ttilde / \beta$) in the the case of covariance mismatch.}
\label{fig:impact_Pfa_Kelly_K=32}
\end{figure}

We now address the main issue of this paper, namely the impact of covariance mismatch on the probability of false alarm of CFAR detectors. For illustration purposes we consider Kelly's detector ($\ttilde$) and the AMF ($\tAMF=\ttilde / \beta$). For both detectors, we set the threshold to ensure $\Pfabar=10^{-3}$ in the absence of mismatch. We generate $50$ different matrices $\mSigmat$ and for each of them we evaluate the actual $\Pfa$ in the presence of covariance mismatch. The results are given in Figure \ref{fig:impact_Pfa_Kelly_K=32}. They confirm that in the presence of covariance mismatch it is actually very complicated to have a constant $\Pfa$ equal to its nominal value. One can observe a significant increase on the mean value of $\Pfa$, especially for the AMF detector which is less robust than Kelly. One can also observe that the increase is more significant in case 1 than in case 2. In fact it seems that the case $\mSigmat = \mSigma^{1/2} \Wt^{-1} \mSigma^{H/2}$ induces much more distortion on $\mSigma$ than when $\mSigma=\U\mLambda\U^{H}$,  $\mSigmat=\U\mLambda^{1/2}\diag{\gamma_{n}}\mLambda^{1/2}\U^{H}$. In the latter case $\mSigma$ and $\mSigmat$ differ only by a``diagonal'' matrix $\diag{\gamma_{n}}$ whereas in the former case the difference is a ``full'' matrix $\Wt$. Whatever the case it appears that the false alarm rate is significantly impacted by covariance mismatch.

\section{Possible mitigation \label{sec:mitigation}}
Given the fluctuations and the increase of $\Pfa$ observed previously,  we now investigate whether it is possible to somewhat mitigate these deleterious effects. As said before, a covariance mismatch has three main consequences, one is the modification of the distribution of $\beta$, the second is the coefficient in $\ttilde$ before the F distribution and the third is the fact that the F distribution becomes non central under $H_{0}$. As for the impact on $\beta$ it is difficult to figure out a remedy, except to find a form that is less sensitive to deviations from the beta distribution. As for the coefficient in \eqref{pdf_ttilde_mismatch}, it is given by $(\Omega_{2.1} + \vxtilde_{1}^{H}\W_{11}^{-1}\vxtilde_{1})/(1+\vxtilde_{1}^{H}\W_{11}^{-1}\vxtilde_{1})$ where $\Omega_{2.1} = (\vv^{H}\mSigmat^{-1}\vv)/(\vv^{H}\mSigma^{-1}\vv)$. The value $\mathtt{1}$ there originates from the fact that $\ttilde = s_{2}/(\mathtt{1}+s_{1}-s_{2})$. This suggests that we could consider
\begin{equation}\label{ttilde_kappa}
\ttilde(\kappa) = \frac{s_{2}}{\kappa+s_{1}-s_{2}} = \frac{|\vx^{H}\St^{-1}\vv|^{2}/(\vv^{H}\St^{-1}\vv)}{\kappa+\vx^{H}\St^{-1}\vx- \frac{|\vx^{H}\St^{-1}\vv|^{2}}{\vv^{H}\St^{-1}\vv}}
\end{equation}
Interestingly enough this is exactly Kalson's detector \cite{Kalson92}. Using Appendix \ref{app:proof} it can be shown that 
\begin{align}\label{pdf_ttilde_kappa_mismatch}
&\ttilde(\kappa) | \vxtilde_{1},\W_{11} = \frac{W_{2.1}^{-1} |\xtilde_{2}-\W_{21}\W_{11}^{-1}\vxtilde_{1}|^{2} }{\kappa+\vxtilde_{1}^{H}\W_{11}^{-1}\vxtilde_{1}} \nonumber \\
&\dist \frac{1+\beta(\Omega_{2.1}-1)}{1+\beta(\kappa-1)} \nonumber \\	
&\, \times \CF{1,K-N+1}{\beta \frac{||\alpha|(\vv^{H}\mSigmat^{-1}\vv)^{1/2}+\mOmega_{21}\mOmega_{11}^{-1}\vxtilde_{1}|^{2}}{1+\beta(\Omega_{2.1}-1)}}
\end{align}
in the general case and
\begin{equation}
\ttilde(\kappa) | \beta \dist \frac{1+\beta(\Omega_{2.1}-1)}{1+\beta(\kappa-1)} \CF{1,K-N+1}{\frac{\beta |\alpha|^{2}(\vv^{H}\mSigmat^{-1}\vv)}{1+\beta(\Omega_{2.1}-1)}}
\end{equation}
in the GER case. Therefore, if $\kappa = \Omega_{2.1} $, the coefficient is equal to $1$ and $\ttilde(\Omega_{2.1})$ has a complex F distribution. Moreover, if the GER holds $\ttilde(\Omega_{2.1}) \dist \CF{1,K-N+1}{0}$ under $H_{0}$ and this detector is CFAR and has the same distribution as $\ttilde$ in the case of no mismatch. Of course the clairvoyant choice $\kappa = \Omega_{2.1}$ cannot be made in practice where one must fix the value of $\kappa$. However the above discussion suggests that a detector of the form \eqref{ttilde_kappa} can somehow mitigate the effect of covariance mismatch. Finally, for the fact that the F distribution becomes non central under $H_{0}$, its consequence is that $\ttilde$ tends to increase in a random (at least unknown) way, and so the probability of false alarm will also increase. To counteract this effect, one could preferably choose $\kappa >  \Omega_{2.1}$ so that a decrease in the coefficient before the F distribution could compensate the increase of the latter term.

\begin{figure*}[htb]
	\centering
	\subfigure{%
		\includegraphics[width=7.5cm]{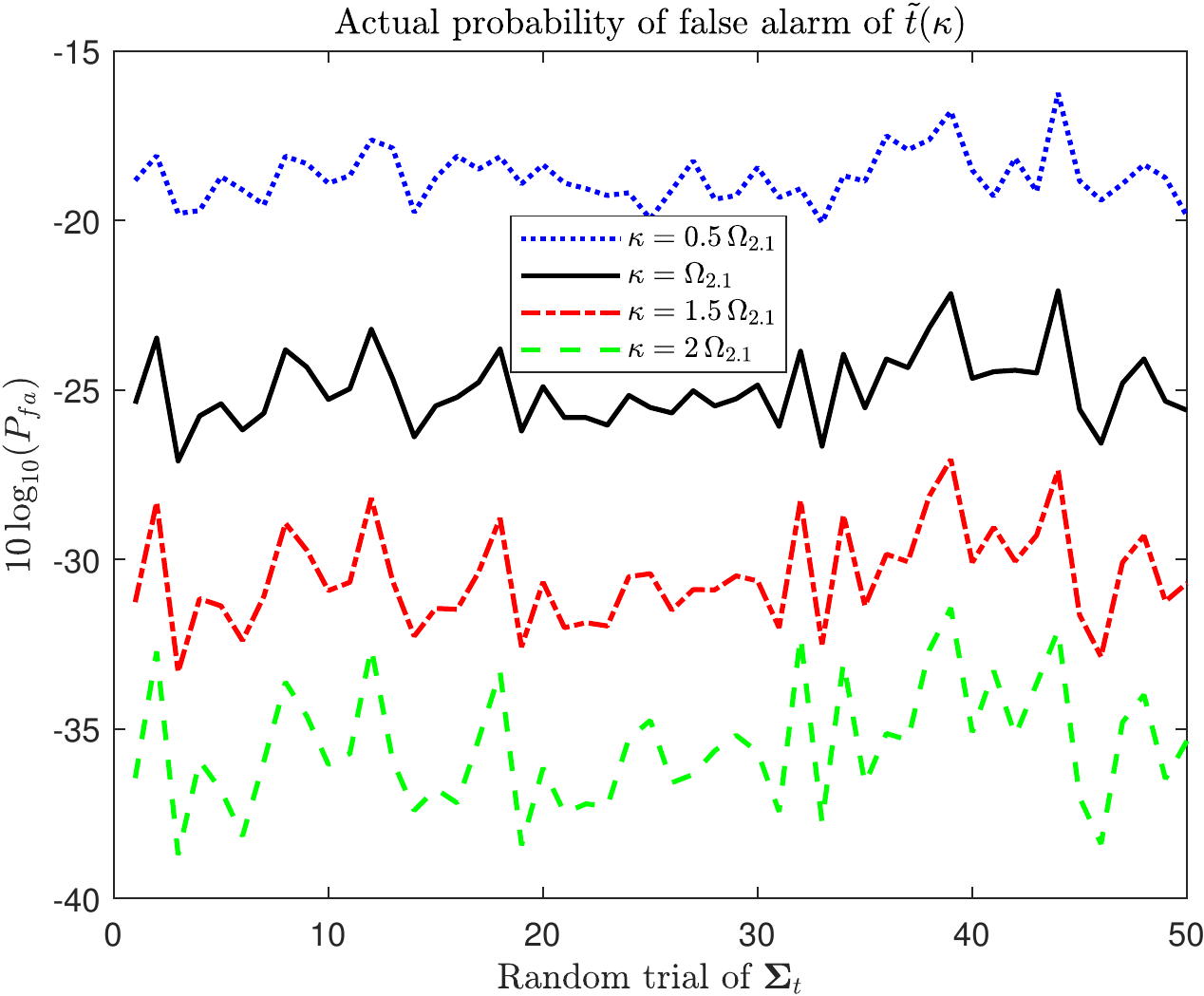}}
	\subfigure{%
		\includegraphics[width=7.5cm]{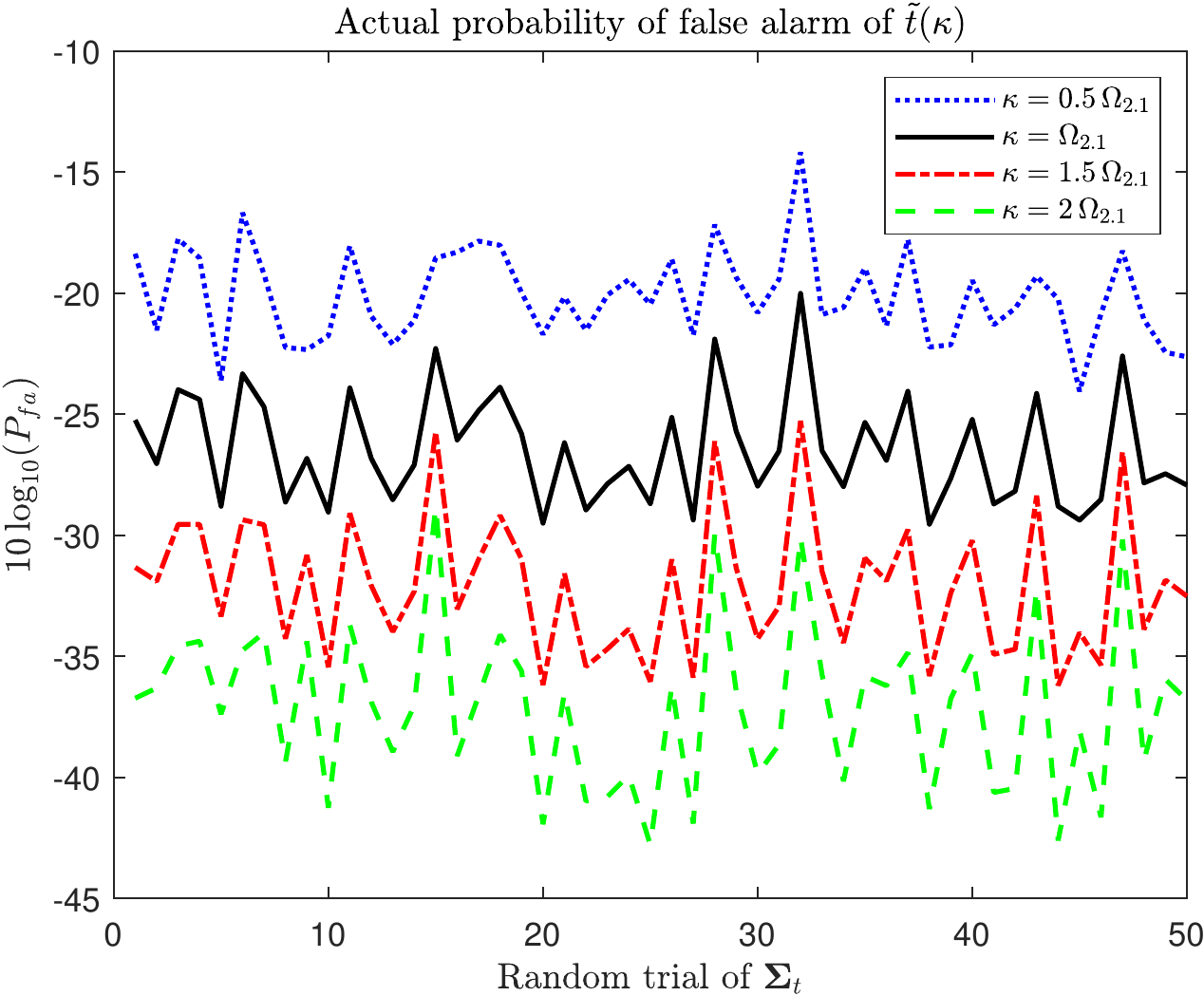}} \\			
	\subfigure{
		\includegraphics[width=7.5cm]{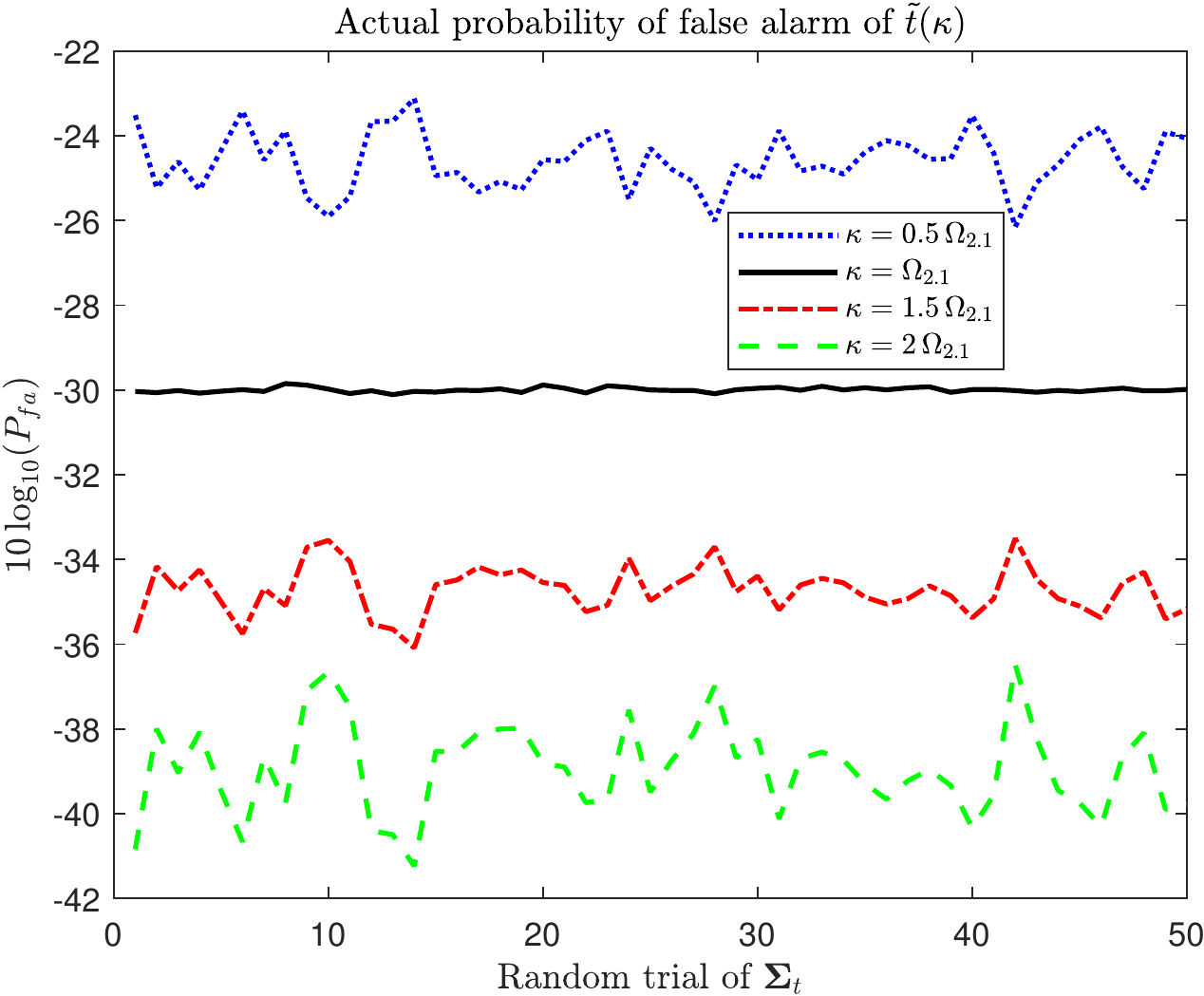}}
	\subfigure{%
		\includegraphics[width=7.5cm]{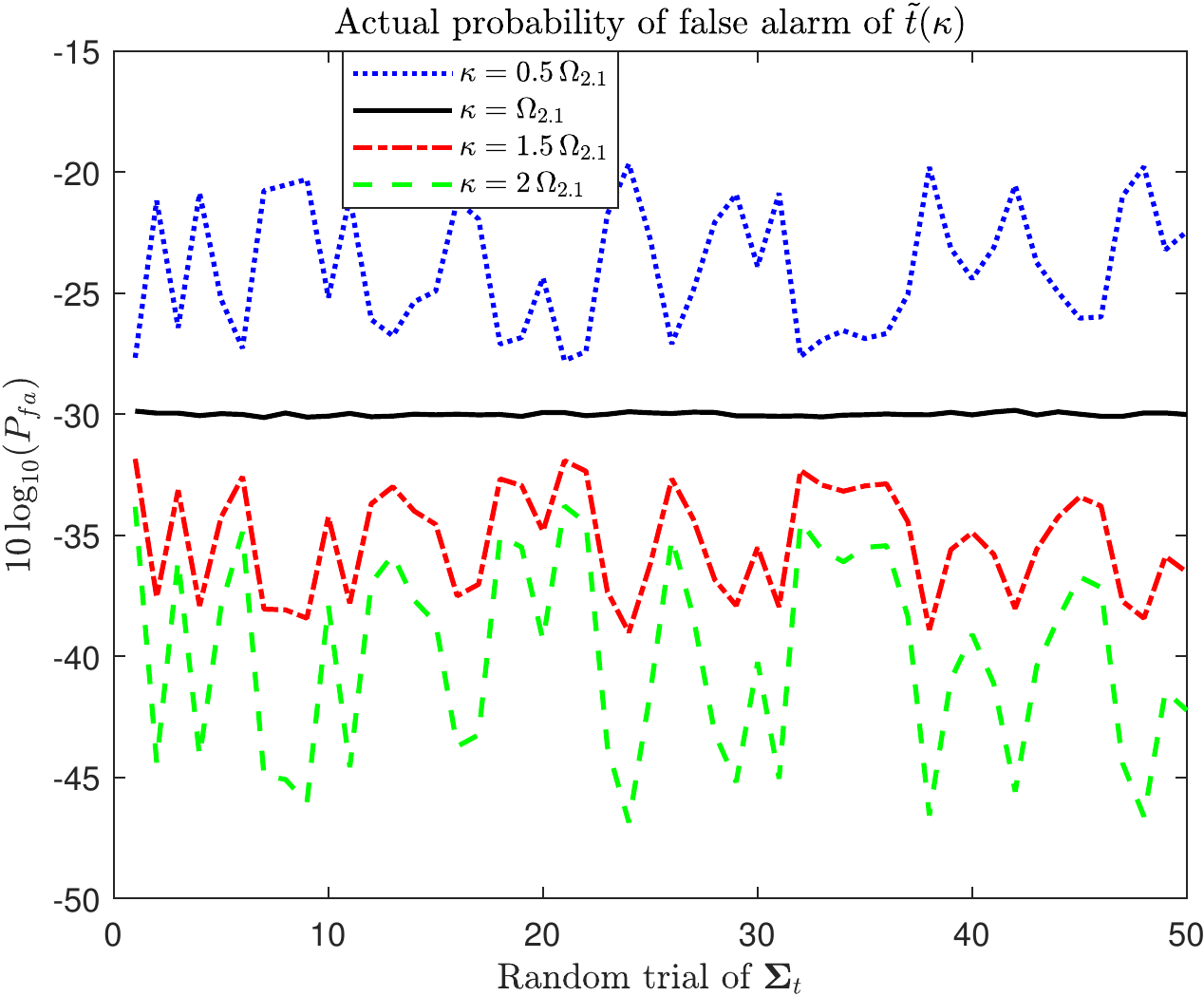}}
	\caption{Actual probability of false alarm of $\ttilde(c \Omega_{2.1})$ in the case of covariance mismatch. The left panel concerns case 1, the right panel case 2. The lower panel corresponds to the case where the GER is satisfied.}
	\label{fig:impact_Pfa_Kalson_adaptive_kappa_K=32}
\end{figure*}

In order to assess the validity of $\ttilde(\kappa)$ we first begin by assuming that $\Omega_{2.1}$ is known and $\kappa$ is selected as $c \Omega_{2.1}$ in order to figure out how precisely the latter should be known. As before, we generated $50$ different matrices $\mSigmat$ and we evaluated $\Prob{\ttilde(c \Omega_{2.1}) \geq \bar{\eta} | \mSigmat \neq \mSigma}$. We want to check that $\Prob{\ttilde(\Omega_{2.1} \geq \bar{\eta} | \mSigmat^{-1}\vv = \lambda \mSigma^{-1}\vv} = \Pfabar$ when the GER is satisfied, and we want to evaluate this probability when the GER is not satisfied. Furthermore, we study how $\Prob{\ttilde(c \Omega_{2.1} \geq \bar{\eta} | \mSigmat \neq \mSigma}$ varies with $c$. The results are displayed in Figure \ref{fig:impact_Pfa_Kalson_adaptive_kappa_K=32}. First, they confirm that, if the GER is satisfied, then $\Prob{\ttilde(\Omega_{2.1}) \geq \bar{\eta} | \mSigmat^{-1}\vv = \lambda \mSigma^{-1}\vv} = \Pfabar$ and hence the detector is CFAR. However, when the GER is no longer satisfied the false alarm rate of $\ttilde(\Omega_{2.1}) $ is above $\Pfabar$, typically around $3 \times 10^{-3}$. Next, as could be anticipated, $\Prob{\ttilde(c \Omega_{2.1}) \geq \bar{\eta} | \mSigmat \neq \mSigma} $ decreases when $c$ increases. Another important feature is that the false alarm rate seems to vary much less when $\mSigmat$ varies, which is a good sign and could indicate that one could reach a close to constant false alarm. To conclude, if $\Omega_{2.1}$ were known, the choice $\kappa=\Omega_{2.1}$ ensures a CFAR if the GER is satisfied and a false alarm rate above its nominal value if the GER is not satisfied. In this case, it is preferable to choose $\Omega_{2.1}  \leq \kappa \leq 2\Omega_{2.1}$ in order to have a probability of false alarm close to $\Pfabar$. 

\begin{figure*}[htb]
	\centering
	\subfigure{%
		\includegraphics[width=7.5cm]{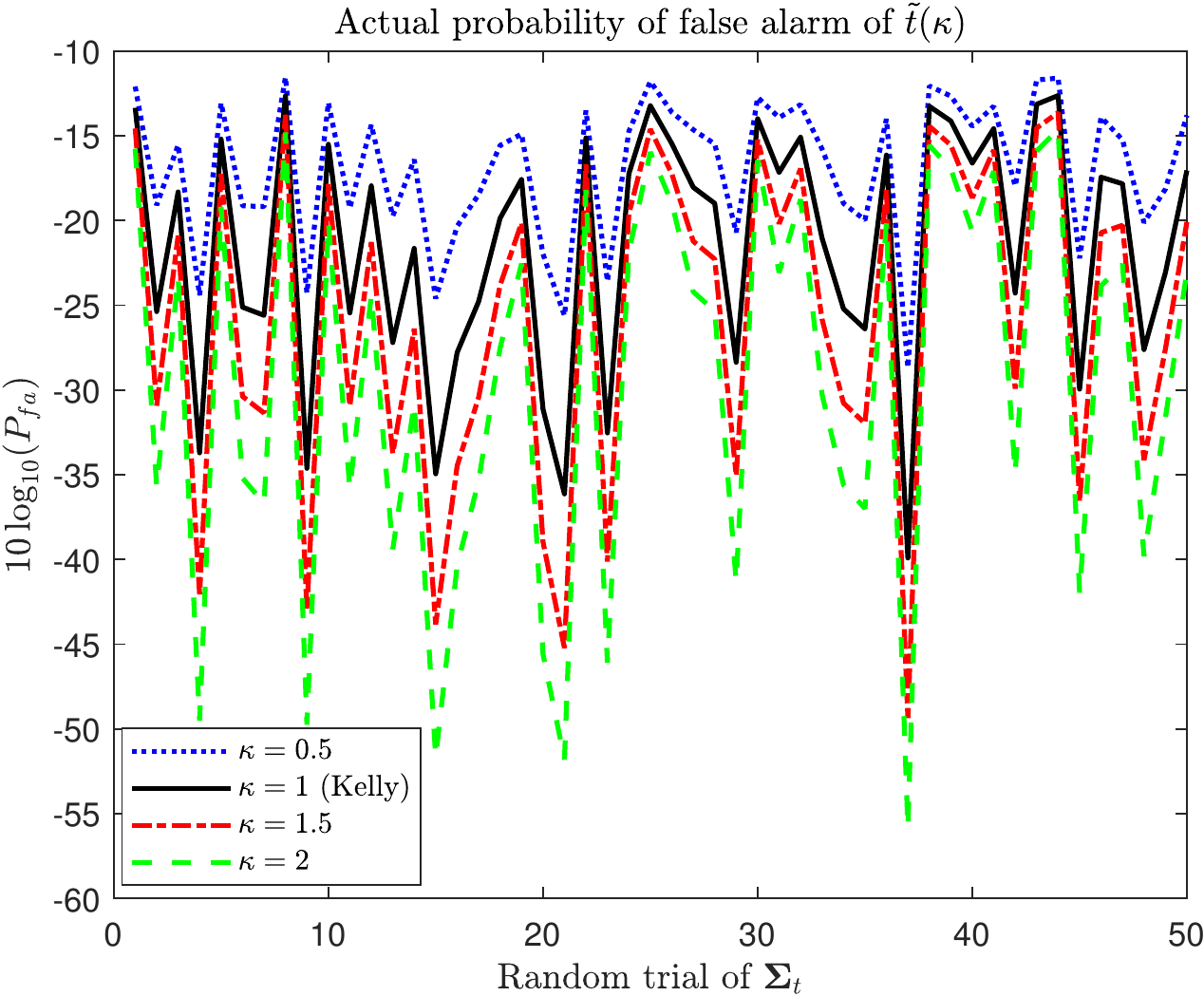}}
	\subfigure{%
		\includegraphics[width=7.5cm]{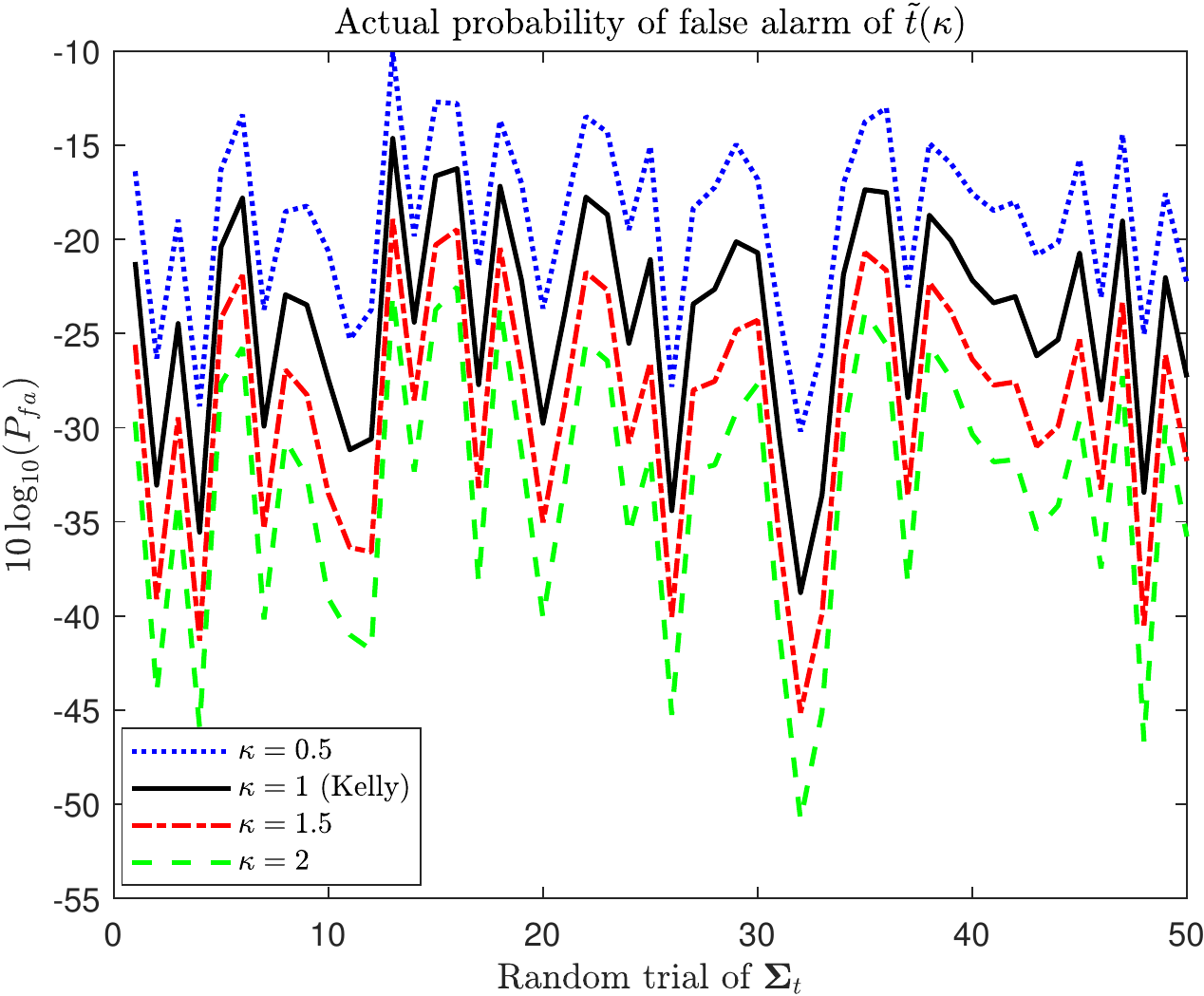}} \\		
	\subfigure{
		\includegraphics[width=7.5cm]{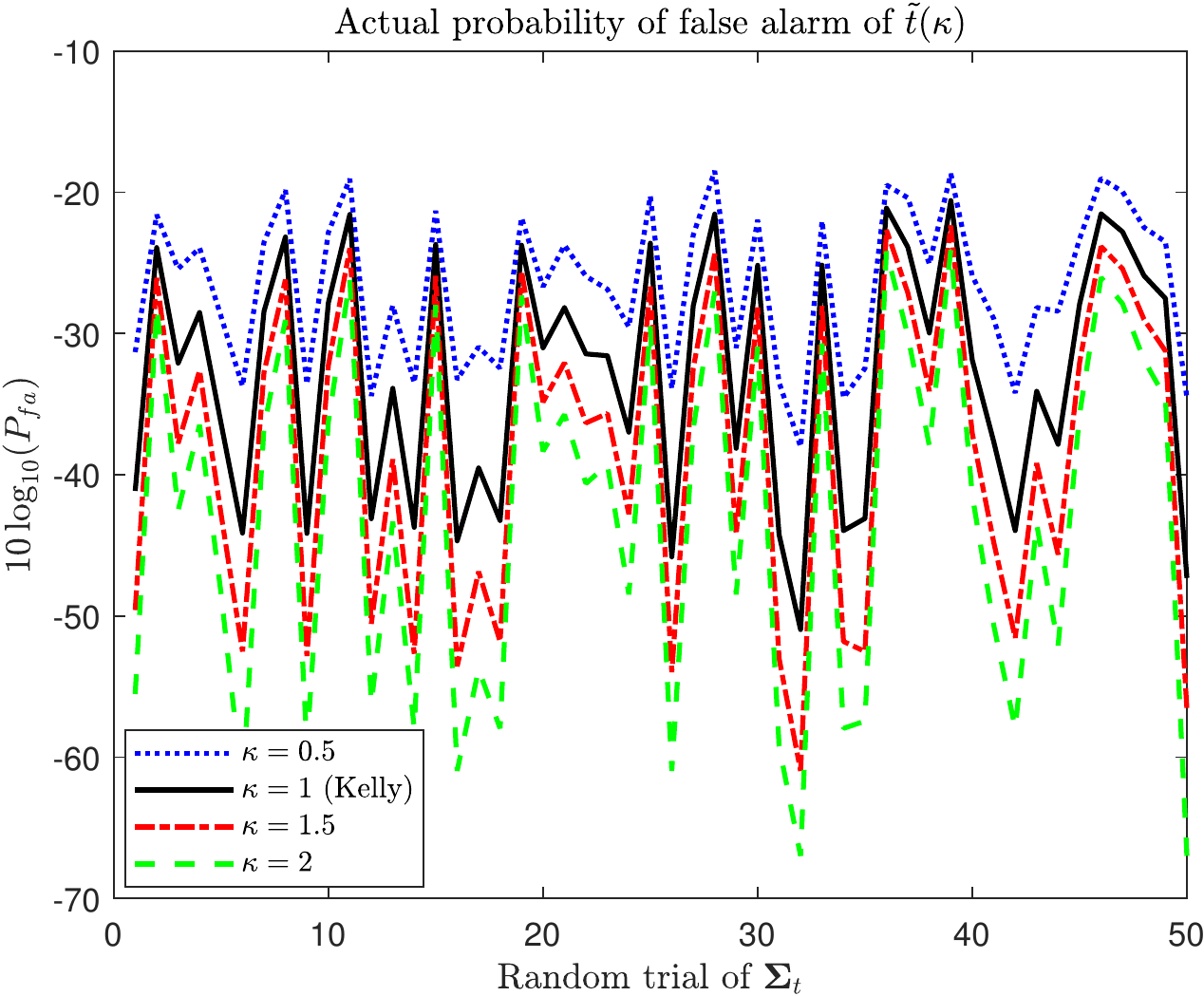}}
	\subfigure{%
		\includegraphics[width=7.5cm]{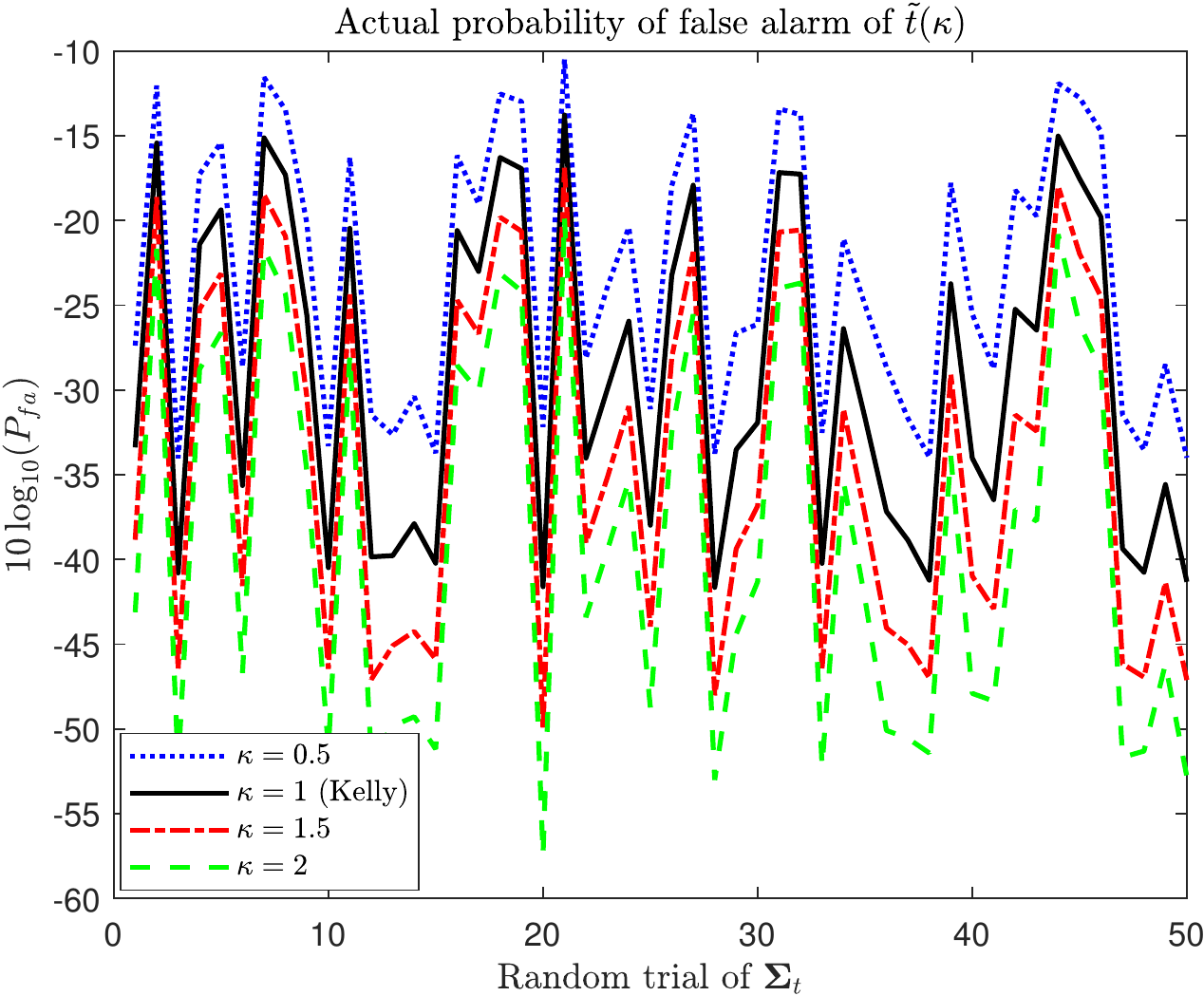}}
	\caption{Actual probability of false alarm of $\ttilde(\kappa)$ in the case of covariance mismatch. The left panel concerns case 1, the right panel case 2. The lower panel corresponds to the case where the GER is satisfied.}
	\label{fig:impact_Pfa_Kalson_fixed_kappa_K=32}
\end{figure*}

\begin{figure}[htb]
\centering
\subfigure[$\Delta=3$dB]{%
\includegraphics[width=8cm]{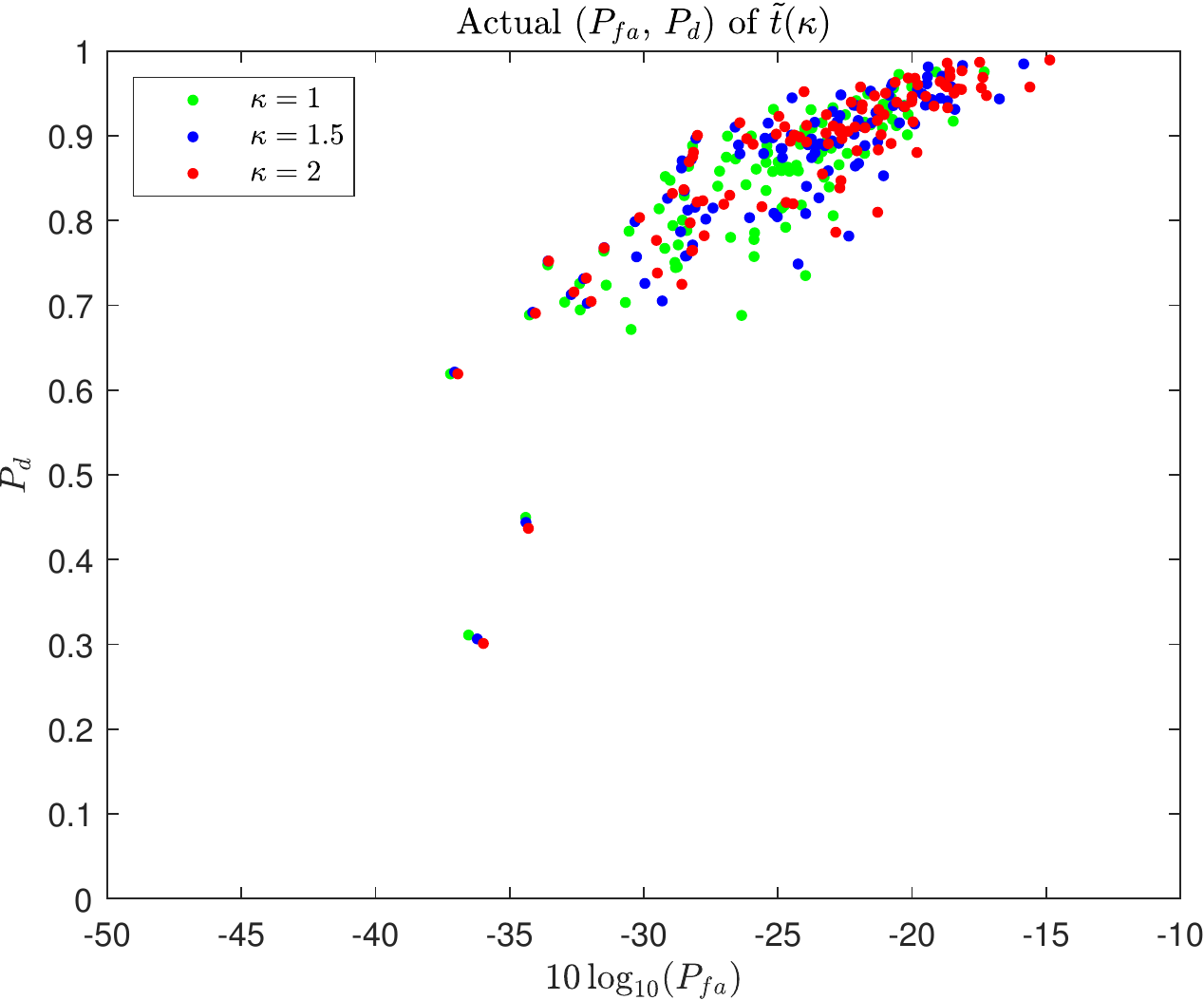}} \\
\subfigure[$\Delta=6$dB]{%
\includegraphics[width=8cm]{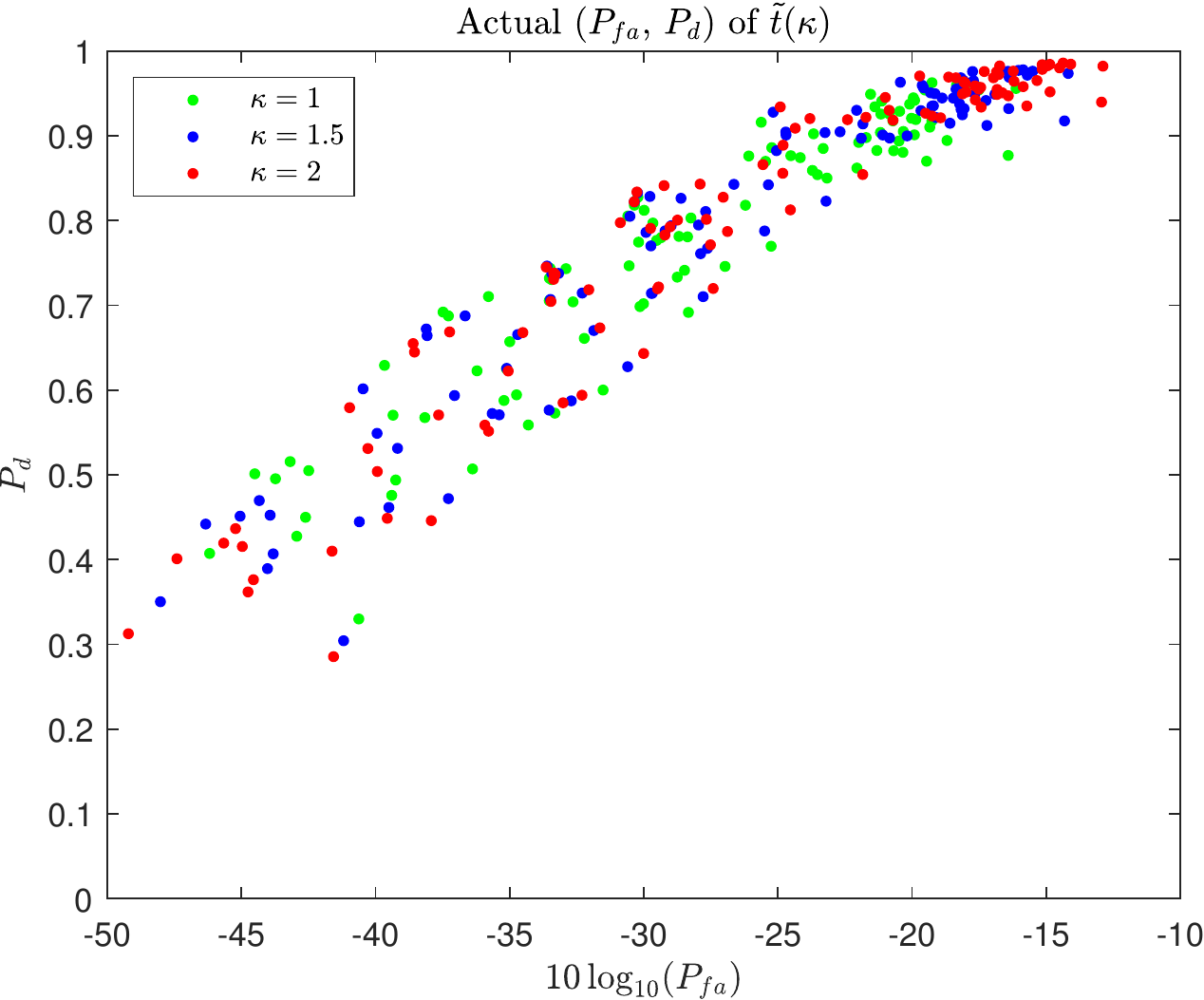}}
\caption{Actual probability of false alarm and probability of detection of $\ttilde(\kappa)$ in case 1. $\mSigmat = \mSigma^{1/2} \Wt^{-1} \mSigma^{H/2}$ with $\E{\Wt^{-1}}=\gamma \I$ 
	and $10\log_{10}\gamma$ uniformly distributed on $[-\Delta,+\Delta]$.}
\label{fig:impact_ROC_Kalson_iW_K=32}
\end{figure}

\begin{figure}[htb]
\centering
\subfigure[$\Delta=3$dB]{%
\includegraphics[width=8cm]{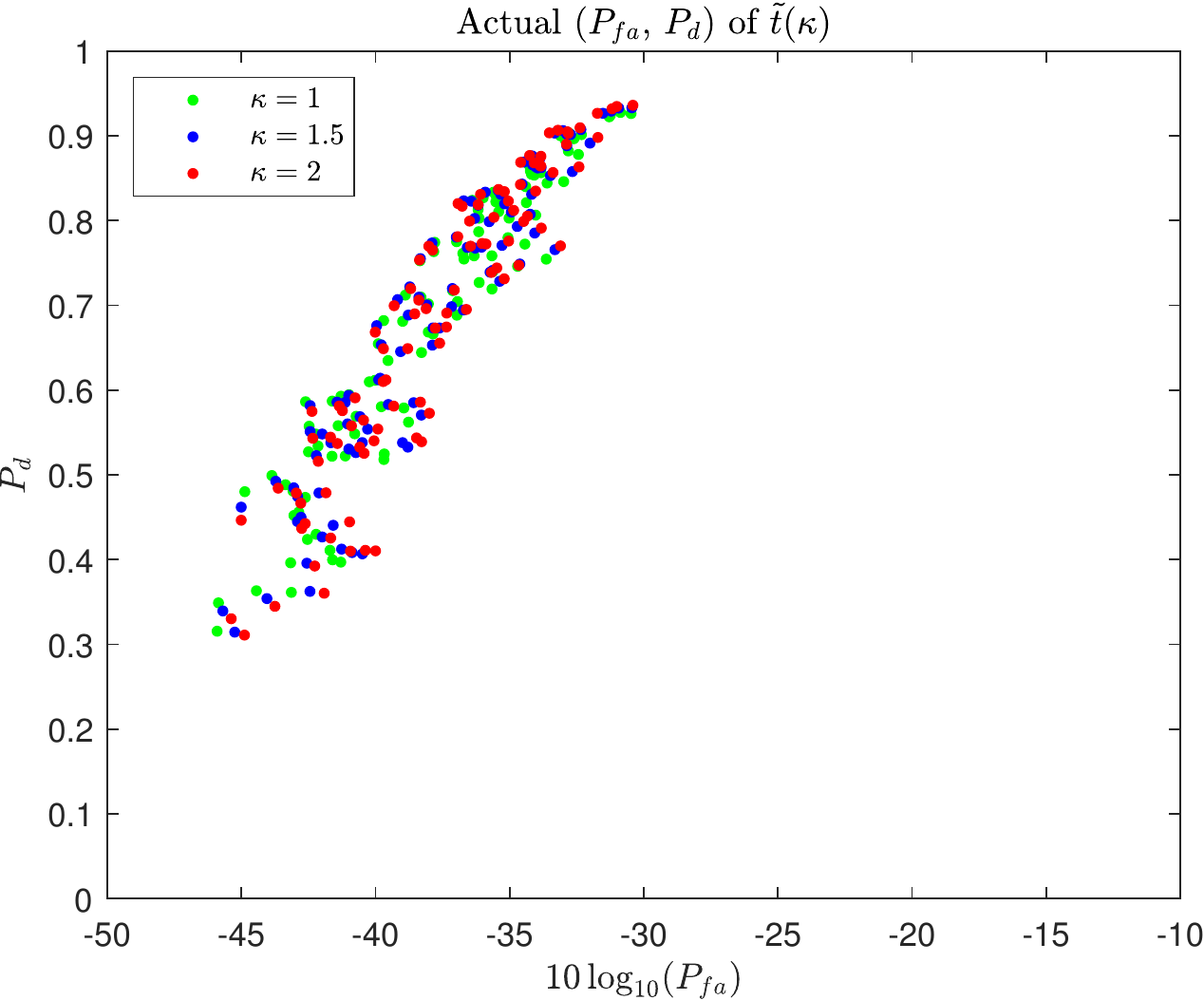}} \\
\subfigure[$\Delta=6$dB]{%
\includegraphics[width=8cm]{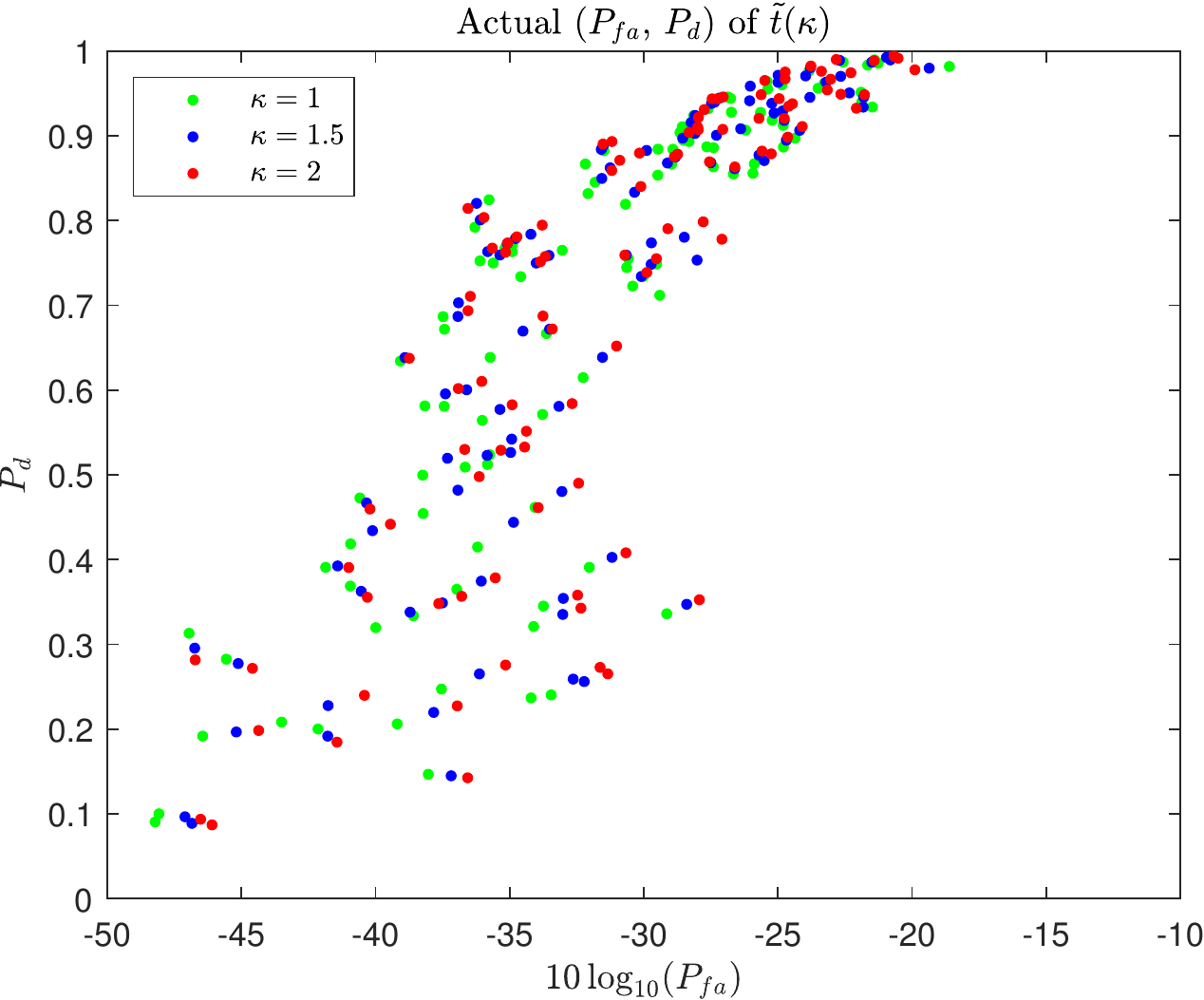}}
\caption{Actual probability of false alarm and probability of detection of $\ttilde(\kappa)$ in case 2. $\mSigma=\U\mLambda\U^{H}$ and $\mSigmat=\U\mLambda^{1/2}\diag{\gamma_{n}}\mLambda^{1/2}\U^{H}$ with $10\log_{10}\gamma_{n}$ uniformly distributed on $[-\Delta,+\Delta]$.}
\label{fig:impact_ROC_Kalson_Lambda_K=32}
\end{figure}

However, in practice $\Omega_{2.1}$ is unknown and $\kappa$ needs to be fixed.  Therefore, in a second simulation we no longer assume that $\Omega_{2.1}$ is known and $\kappa$ is set to a fixed value. Again we evaluate $\Prob{\ttilde(\kappa) \geq \bar{\eta} | \mSigmat \neq \mSigma}$ for $50$ different matrices $\mSigmat$. Figure \ref{fig:impact_Pfa_Kalson_fixed_kappa_K=32} clearly indicates that a fixed value of $\kappa$ does not completely fix the problem. Indeed, one can observe significant variations of $\Prob{\ttilde(\kappa) \geq \bar{\eta} | \mSigmat \neq \mSigma}$. However, with $\kappa=2$, we get an average value of this probability close to $10^{-3}$ with a standard deviation of about $10^{-3}$. This means that the choice $\kappa=2$ provides a detector that is less sensitive to the fact that $\mSigmat \neq \mSigma$ and to different outcomes of $\mSigmat$ and thus is helpful in controlling $\Pfa$. With this choice, the actual false alarm rate does not deviate much from its target value.

Finally, while we concentrated so far about $\Pfa$, we now investigate the impact of covariance mismatch on $\Pd(\ttilde(\kappa) | \mSigmat \neq \mSigma)$. As before, the target $\Pfa$ is $10^{-3}$ and, to anticipate the increase of $\Pfa$ due to covariance mismatch, the threshold $\eta(\kappa)$ is set so that $\Pfabar = \Proba{H_0}{\ttilde(\kappa) \geq \eta(\kappa) | \mSigmat = \mSigma} = 10^{-4}$. Similarly, we set the SNR so that $\Pdbar = \Proba{H_1}{\ttilde(\kappa) \geq \eta(\kappa) | \mSigmat = \mSigma}= 0.7$. Our aim is to evaluate the fluctuations of $\Pfa(\kappa)=\Proba{H_0}{\ttilde(\kappa) \geq \eta(\kappa) | \mSigmat \neq \mSigma}$ and $\Pd(\kappa) = \Proba{H_1}{\ttilde(\kappa) \geq \eta(\kappa) | \mSigmat \neq \mSigma}$ around the operating point $(\Pfabar,\Pdbar)$ of the ROC. The results are displayed in Figures \ref{fig:impact_ROC_Kalson_iW_K=32}-\ref{fig:impact_ROC_Kalson_Lambda_K=32} for two different levels of mismatch. It can be observed that $\Pfa$ is more impacted as we most often have $\Pfa > \Pfabar$ and the variations may be important, up to 2-3 decades. On the other hand, $\Pd$ is less impacted: for instance when the actual $\Pfa \simeq 10^{-4}$ then $\Pd \simeq 0.5-0.6$. And, as already said before, $\Pd$ varies  simply because $\Pfa$ varies a lot, which confirms that control of the false alarm rate is primordial.

\section{Conclusion \label{sec:conclusion}}
In this paper we considered a possible covariance mismatch between the data under test and the training samples and we investigated its consequences on the probability of false alarm of detectors which depend on the variables $\beta$ and $\ttilde$ and hence are CFAR in the case of no mismatch. A statistical representation of these two variables was obtained for an arbitrary covariance mismatch and for both the null and alternative hypothesis. It was shown that $\beta$ is no longer beta distributed and that $\ttilde$ is strongly impacted. Numerical simulations illustrated the important variation of the probability of false alarm when $\mSigmat$ varies. We showed that, when the GER is satisfied, a clairvoyant modification of Kelly's detector yields a CFAR detector. We investigated Kalson's detector to mitigate covariance mismatch effects and showed that setting $\kappa=1.5$, $\kappa=2$ seems a good choice to mitigate variations of the false alarm rate. Yet, it appears nearly impossible to obtain a strict constant false alarm rate when $\mSigmat \neq \mSigma$ with the current data model. A possible solution lies in assuming that more information about $\mSigma$ could be available. This is the approach taken in \cite{Raghavan19b}  where  an additional set of training samples with covariance matrix $\mSigma$ is available, and where a CFAR detector based on cell-averaging is derived. In \cite{Besson20d} we derived a detector which has better detection than cell averaging but is not strictly CFAR. Therefore, a perspective is to identify statistics of this new model whose distribution is parameter free under $H_{0}$.

\appendix
\section{Proof of \eqref{storep_beta_mismatch}-\eqref{storep_ttilde_mismatch} \label{app:proof}}
In this appendix we derive the stochastic representation of $(\beta,\ttilde)$. Let $\Vorth$ be a semi-unitary matrix orthogonal to $\vv$, i.e., $\Vorth^{H}\Vorth =\eye{N-1}$ and $\Vorth^{H}\vv=\vzero$ and let us assume without loss of generality that $\left\| \vv  \right\|=1$. Let $\mSigma^{1/2}$ and $\mSigmat^{1/2}$ denote square-roots of $\mSigma$ and $\mSigmat$. First note that $\vx \dist  (\alpha \vv + \mSigma^{1/2}\vn)$ with $\vn \dist \vCN{N}{\vzero}{\eye{N}}$ and that  $\St \dist  \mSigmat^{1/2} \Wt \mSigmat^{H/2}$ where $\Wt \dist  \CW{N}{K}{\eye{N}} $ follows a complex Wishart distribution with $K$ degrees of freedom and parameter matrix $\eye{N}$. 
First, let us rewrite $s_{1}$ as
\begin{align}
	s_{1} & =\vx^{H}\St^{-1}\vx \nonumber \\
	&= \vx^{H} \mSigmat^{-H/2} \Wt^{-1} \mSigmat^{-1/2}  \vx \nonumber \\
	&=  (\Q^{H} \mSigmat^{-1/2}\vx)^{H}  \Q^{H}\Wt^{-1}\Q (\Q^{H} \mSigmat^{-1/2}\vx)
\end{align}
where $\Q$ is a non-singular matrix. Let $\Gt$ be a square-root of $\mSigmat$, i.e., $\mSigmat=\Gt\Gt^{H}$ and let us choose $\Q = \begin{bmatrix} \Gt^{H} \Vorth \Ft^{-H} & (\vv^{H}\mSigmat^{-1}\vv)^{1/2} \Gt^{-1} \vv \end{bmatrix}$ where $\Ft = (\Vorth^{H} \mSigmat \Vorth)^{1/2}$  so that $\Q$ is unitary and $\Q^{H}\mSigmat^{-1/2}\vv = (\vv^{H}\mSigmat^{-1}\vv)^{1/2}\elast$ with $\elast = \begin{bmatrix} \vzero^{T} & 1 \end{bmatrix}^{T}$. Then, noting that $\W=\Q^{H}\Wt\Q \dist \CW{N}{K}{\eye{N}} $, we obtain
\begin{align}
	s_{1} &\dist \vxtilde^{H}  \W^{-1}  \vxtilde \nonumber \\
	\vxtilde &= \Q^{H} \mSigmat^{-1/2}\vx = (\gammat^{1/2}\elast + \mOmega^{1/2}\vn)
\end{align}
with $\gammat = |\alpha|^{2} (\vv^{H}\mSigmat^{-1}\vv)$ and where 
\begin{align}
	\mOmega = \Q^{H}\mSigmat^{-1/2} \mSigma \mSigmat^{-H/2} \Q
\end{align}
Similarly
\begin{align}
	s_{2} &= \frac{|\vx^{H}\St^{-1}\vv|^{2}}{\vv^{H}\St^{-1}\vv} \nonumber \\
	&= \frac{| (\Q^{H} \mSigmat^{-1/2}\vx)^{H}  \Q^{H}\Wt^{-1}\Q (\Q^{H} \mSigmat^{-1/2}\vv)|^{2}}{ (\Q^{H} \mSigmat^{-1/2}\vv)^{H}  \Q^{H}\Wt^{-1}\Q (\Q^{H} \mSigmat^{-1/2}\vv)} \nonumber \\
	&= \frac{|\vxtilde^{H}\W^{-1}\elast|^{2}}{\elast^{H}\W^{-1}\elast}
\end{align}
Let us partition $\vxtilde$ as $\vxtilde = \begin{pmatrix} \vxtilde_{1} \\ \xtilde_{2}  \end{pmatrix} $ and $\W$ as 
\begin{equation}\label{partition_W}
	\W = \begin{pmatrix} \W_{11} & \W_{12} \\ \W_{21} & W_{22} \end{pmatrix}
\end{equation}
where $ \W_{11}$ is $(N-1)\times(N-1)$. Using the fact that
\begin{equation}
	\W^{-1} = W_{2.1}^{-1} \begin{pmatrix} -\W_{11}^{-1}\W_{12} \\ 1 \end{pmatrix}  \begin{pmatrix} -\W_{21}\W_{11}^{-1} & 1 \end{pmatrix}  + \begin{pmatrix} \W_{11}^{-1} & \vzero \\ \vzero & 0 \end{pmatrix}
\end{equation}
with $W_{2.1}=\W_{22} - \W_{21}\W_{11}^{-1} \W_{12 }$, we can write
\begin{align}\label{representation_s1}
	s_{1} &\dist \vxtilde^{H}  \W^{-1}  \vxtilde \nonumber \\
	&= W_{2.1}^{-1} |\xtilde_{2}-\W_{21}\W_{11}^{-1}\vxtilde_{1}|^{2} + \vxtilde_{1}^{H}\W_{11}^{-1}\vxtilde_{1}
\end{align}
and
\begin{align}\label{representation_s2}
	s_{2} &= \frac{|\vxtilde^{H}\W^{-1}\elast|^{2}}{\elast^{H}\W^{-1}\elast} \nonumber \\
	&= W_{2.1}^{-1} |\xtilde_{2}-\W_{21}\W_{11}^{-1}\vxtilde_{1}|^{2} 
\end{align}
Since $\W_{11} \dist \CW{N-1}{K}{\eye{N-1}}$, it follows that
\begin{equation}
	\vxtilde_{1}^{H}\W_{11}^{-1}\vxtilde_{1} \dist \frac{\vxtilde_{1}^{H}\vxtilde_{1}}{\Cchisquare{K-N+2}{0}}
\end{equation}
Now, $\vxtilde = \Q^{H} \mSigmat^{-1/2}\vx \dist \vCN{N}{\gammat^{1/2}\elast}{\mOmega}$, which implies that $\vxtilde_{1} \dist \vCN{N-1}{\vzero}{\mOmega_{11}}$ and hence $\vxtilde_{1} \dist \mOmega_{11}^{1/2} \vn_{1}$ with $\vn_{1} \dist \vCN{N-1}{\vzero}{\eye{N-1}}$ .  Therefore, a stochastic representation for $\beta$ can be obtained as
\begin{align}\label{pdf_beta_mismatch}
	\beta  &= (1+\vxtilde_{1}^{H}\W_{11}^{-1}\vxtilde_{1})^{-1} \nonumber \\
	&\dist \left[1+\frac{\vn_{1}^{H}\mOmega_{11}\vn_{1}}{\Cchisquare{K-N+2}{0}}\right]^{-1} \nonumber \\
	&\dist \left[1+\frac{\sum_{i=1}^{N-1}\lambda_{i} \Cchisquare{1}{0}}{\Cchisquare{K-N+2}{0}}\right]^{-1}
\end{align}
where $\lambda_{i}$ ($i=1,\ldots,N-1$) are the eigenvalues of $\mOmega_{11}$.

Let us now turn to analysis of $\ttilde$. Since $\W \dist \CW{N}{K}{\eye{N}} $ it can be written as $\W = \Z \Z^{H}$ where $\Z \dist \mCN{N}{\Mzero}{\eye{N}}{\eye{K}}$. Let us partition $\Z = \begin{pmatrix} \Z_{1} \\ \Z_{2} \end{pmatrix}$ so that $\W_{11} = \Z_{1}\Z_{1}^{H}$ and $\W_{21} = \Z_{2}\Z_{1}^{H}$. It ensues that
\begin{equation}\label{key}
	\xtilde_{2}-\W_{21}\W_{11}^{-1}\vxtilde_{1} = \xtilde_{2}-\Z_{2}\Z_{1}^{H}(\Z_{1}\Z_{1}^{H})^{-1}\vxtilde_{1}
\end{equation}
Since $\vxtilde \dist \vCN{N}{\gammat^{1/2}\elast}{\mOmega}$, one has
\begin{align}
	&\xtilde_{2} | \vxtilde_{1} \dist \CN{\gammat^{1/2}+\mOmega_{21}\mOmega_{11}^{-1}\vxtilde_{1}}{\Omega_{2.1}} \\
	&\Z_{2}\Z_{1}^{H}(\Z_{1}\Z_{1}^{H})^{-1}\vxtilde_{1} | \vxtilde_{1},\W_{11}  \dist \CN{0}{\vxtilde_{1}^{H}\W_{11}^{-1}\vxtilde_{1}}
\end{align}
and hence
\begin{align}
	&\xtilde_{2}-\W_{21}\W_{11}^{-1}\vxtilde_{1} | \vxtilde_{1},\W_{11} \nonumber \\
	&\dist \CN{\gammat^{1/2}+\mOmega_{21}\mOmega_{11}^{-1}\vxtilde_{1}}{\Omega_{2.1} + \vxtilde_{1}^{H}\W_{11}^{-1}\vxtilde_{1}} 
\end{align}
Noting that $\W_{2.1} \dist \Cchisquare{K-N+1}{0}$, it follows that
\begin{align}\label{pdf_ttilde_mismatch}
	\ttilde | \vxtilde_{1},\W_{11} &= \frac{W_{2.1}^{-1} |\xtilde_{2}-\W_{21}\W_{11}^{-1}\vxtilde_{1}|^{2} }{1+\vxtilde_{1}^{H}\W_{11}^{-1}\vxtilde_{1}} \nonumber \\
	&\dist \frac{\Omega_{2.1} + \vxtilde_{1}\W_{11}^{-1}\vxtilde_{1}}{1+\vxtilde_{1}^{H}\W_{11}^{-1}\vxtilde_{1}} \nonumber \\
	&\, \times	\CF{1,K-N+1}{\frac{|\gammat^{1/2}+\mOmega_{21}\mOmega_{11}^{-1}\vxtilde_{1}|^{2}}{\Omega_{2.1} + \vxtilde_{1}^{H}\W_{11}^{-1}\vxtilde_{1}}}
\end{align}
Using the fact that $\Omega_{2.1} = (\vv^{H}\mSigmat^{-1}\vv)/(\vv^{H}\mSigma^{-1}\vv)$ and $\vxtilde_{1}^{H}\W_{11}^{-1}\vxtilde_{1} = (1-\beta)/\beta$, it is straightforward to show that \eqref{pdf_ttilde_mismatch} can be rewritten as \eqref{storep_ttilde_mismatch}, which concludes the proof.

\section{The GER for $\mSigma=\U\mLambda\U^{H}$ and $\mSigmat=\U\mLambda^{1/2}\Wt\mLambda^{1/2}\U^{H}$\label{app:GER}}
In section \ref{sec:analysis}, we considered a mismatch between the eigenvalues of $\mSigma$ and $\mSigmat$ while their eigenvectors are the same, i.e., we had $\mSigma=\U\mLambda\U^{H}$ and $\mSigmat=\U\mLambda^{1/2}\diag{ \gamma_{n}}\mLambda^{1/2}\U^{H}$. Proceeding this way the GER is not satisfied. Hence, in order to enforce the GER and to have an equivalent form we consider $\mSigmat=\U\mLambda^{1/2}\Wt\mLambda^{1/2}\U^{H}$ and we look for $\Wt$ that makes $\mSigmat^{-1}\vv = \ell_{2}\mSigma^{-1}\vv$. One has
\begin{align}
\mSigmat^{-1}\vv = \ell_{2}\mSigma^{-1}\vv &\Leftrightarrow \U\mLambda^{-1/2}\Wt^{-1}\mLambda^{-1/2}\U^{H}\vv = \ell_{2}\U\mLambda^{-1}\U^{H}\vv \nonumber \\
&\Leftrightarrow \Wt^{-1}\mLambda^{-1/2}\U^{H}\vv = \ell_{2}\mLambda^{-1/2}\U^{H}\vv 
\end{align}
which implies that $\mLambda^{-1/2}\U^{H}\vv$ is an eigenvector of $\Wt$ associated with eigenvalue $\ell_{2}^{-1}$. Therefore the eigenvalue decomposition of $\Wt$ is given by
\begin{align}\label{gen_GER_wLambda}
\Wt &= \V \begin{pmatrix} \vl_{1}^{-1} & \vzero \\ \vzero & \ell_{2}^{-1} \end{pmatrix} \V^{H} \\
\V &= \begin{bmatrix} \mLambda^{1/2}\U^{H}\Vorth\F^{-H} & \frac{\mLambda^{-1/2}\U^{H}\vv}{(\vv^{H}\mSigma^{-1}\vv)^{1/2}} \end{bmatrix}
\end{align}
where $\F=(\Vorth^{H}\mSigma\Vorth)^{1/2}$. Equation \eqref{gen_GER_wLambda} enables one to generate a matrix $\mSigmat$ such that the GER is satisfied. When all eigenvalues are equal to $1$, $\mSigma=\mSigmat$ so the values of $\vl_{1}$ and $\ell_{2}$ determine the degree of mismatch.

\end{document}